%% file: Rpsi.tex
\documentclass[12pt, a4paper,english]{smfart}

\usepackage{pstricks-add}
\usepackage{amsmath}

\input{macros.sty}

\newcommand\myatop[2]{\genfrac{}{}{0pt}{}{#1\hfill}{#2\hfill}}

\theoremstyle{plain}

\title[Nearby cycles in the semi stable case]{Filtrations of the perverse sheaf of nearby cycles in the semi stable situation}

\author{Pascal Boyer}
\email{boyer@math.univ-paris13.fr}
\address{Universit\'e Paris 13, Sorbonne Paris Nord \\
LAGA, CNRS, UMR 7539\\ 
F-93430, Villetaneuse (France) \\
PPAL: ANR-25-PRC}

\begin{document}

\setcounter{tocdepth}{3}

\maketitle

\begin{abstract}
In the strict semi stable reduction situation, we describe the various filtrations of the
perverse sheaf of nearby cycles in terms of irreducible perverse sheaves together
with the action of the monodromy operator. We then study
the spectral sequences associated to these filtration computing the sheaf cohomology groups.
Finally we propose an illustration of how it can be used to compute the cohomology groups.
Considering the similarity with the
results of \cite{boyer-duke}, it could be a good introduction before reading loc. cit.
\end{abstract}

\maketitle

\tableofcontents

\section{Introduction}

Throughout the paper $\Lambda$ will denote either $\overline \Qm_l$, $\overline \Zm_l$
or $\Fm_l$ where $l$ is some prime number distinct from another prime number $p$ equal to the
residue characteristic $\kappa$ of some discrete henselian valuation ring $R$.

The computation of nearby cycles $R^i \Psi(\Lambda)$ for a scheme
$X \longrightarrow \spec R$ with semi-stable reduction, is well understood, 
see \cite{ill} theorem 3.2 for example. 
It is also well known that the complex $\Psi(\Lambda):=R\Psi(\Lambda)[d-1](\frac{d-1}{2})$ of nearby cycles, where $d=\dim X$, is
perverse and autodual relative to the Grothendieck-Verdier duality.
The first goal of this work is  to describe, when $\Lambda=\overline \Qm_l$ or
$\overline \Fm_l$, its irreducible constituants in the category of perverse sheaves on 
$Y$ the geometric special fiber of $X$.

In \cite{boyer-torsion} we explain, for $\Lambda=\overline \Qm_l, \overline \Zm_l$ or
$\overline \Fm_l$, given a stratification of $Y$, how to construct
various filtrations of $\Psi(\Lambda)$. The second goal of this work is then to describe these
filtrations and understand which extensions are split or not.
We finally explain how one can define the nilpotent monodromy operator from our computations.

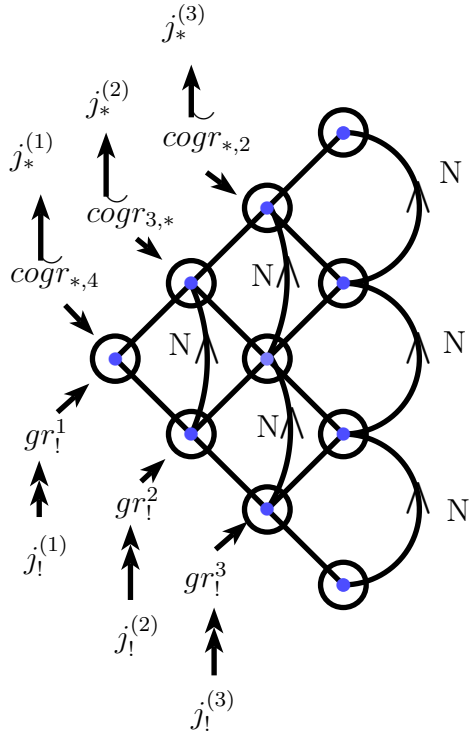
\begin{figure}[htbp]
\begin{center}
\input{fig-ss-N.tex}
\end{center}
\caption{\label{fig-1} Filtrations of $\Psi(\Lambda)$ when $r=4$.}
\end{figure}

The results are illustrated in figure \ref{fig-1} as explained after and using the notations of the next sections
\begin{itemize}
\item Each circle represents a graded part isomorphic to some $i^{(h)}_*\lexp p j^{(h)}_{!*} \Lambda^{(h)} (\frac{h-1-2(k-1)}{2})$ for $1 \leq h \leq r$ and $1 \leq k \leq h$, cf. notation \ref{nota-lambdaI}.

\item The first diagonal line is the first graded part $\gr^1_!(\Psi(\Lambda))$ 
of the filtration $\Fil^\bullet_!(\Psi(\Lambda))$ of stratification of $\Psi(\Lambda)$, cf.
definition \ref{defi-fil!}, constructed
using the adjunction morphisms $j^{(h)}_! j^{(h)} \rightarrow \Id$.
It also admits a filtration $\Fill^{-k}(\gr^1_!(\Psi(\Lambda))$  for $1 \leq k \leq r$,
where the last graded part $\grr_{-1}(\gr^1_!(\Psi(\Lambda))$, i.e. the top quotient, 
corresponds to the circle on the left
of this line and corresponds to $i^{(1)}_* \lexp p j^{(1)}_{!*} \Lambda^{(1)}$, the previous one corresponds
to $i^{(2)}_* \lexp p j^{(2)}_{!*} \Lambda^{(2)} (\frac{1}{2})$ and so on until the first one 
$\grr^{-r}(\gr^1_!(\Psi(\Lambda))$, the socle subspace, is isomorphic to
$i^{(r)}_* \lexp p j^{(r)}_{!*} \Lambda^{(r)}(\frac{r-1}{2})$, cf. \S \ref{para-fil!}.

\item The second diagonal represents $\gr^2_!(\Psi(\Lambda))$ having a filtration 
$\Fill^{-k}(\gr^2_!(\Psi(\Lambda))$ for $2 \leq k \leq r$,
where the last graded parts $\grr^2(\gr^2_!(\Psi(\Lambda))$
 appears on the left of this line and corresponds to
$i^{(2)}_* \lexp p j^{(2)}_{!*} \Lambda^{(2)} (\frac{-1}{2})$, the previous one corresponding to
$i^{(3)}_* \lexp p j^{(3)}_{!*} \Lambda^{(3)}$ and so on until the first one 
$\grr^{-r}(\gr^2_!(\Psi(\Lambda)) \simeq i^{(r)}_* \lexp p j^{(r)}_{!*} \Lambda^{(r)}(\frac{r-3}{2})$.

\item The graded parts have to be read from the top to the bottom, in particular the socle (resp. top)
appears in the top (resp. bottom) of the figure and corresponds to 
$i^{(r)}_* \lexp p j^{(r)} \Lambda^{(r)} (\frac{r-1}{2})$ (resp. 
$i^{(r)}_* \lexp p j^{(r)} \Lambda^{(r)} (-\frac{r-1}{2})$).

\item The nilpotent monodromy operator, defined whatever are the coefficients ring, 
is represented by the arcs of circle starting from bottom to top, cf. \S \ref{para-N}

\item Each segment represents some non trivial extensions in the sense that for a segment with slope
$1$ (resp. $-1$) and for $I \subset J \subset \{ 1,\cdots,r \}$ with $\sharp J \setminus I=1$,
$i_{J,*} \lexp p j_{J,!*} \Lambda_J(\frac{\delta}{2})$ (resp. 
$i_{I,*} \lexp p j_{I,!*} \Lambda_I(\frac{\delta}{2})$)
being an irreducible constituant of the circle in the top of this segment and
$i_{I,*} \lexp p j_{I,!*} \Lambda_I(\frac{\delta-1}{2})$ (resp. 
$i_{J,*} \lexp p j_{J,!*} \Lambda_J(\frac{\delta-1}{2})$)
being an irreducible constituant of the circle of the bottom of this segment,
then the first have to appears before the second one as a graded parts of any filtration of 
$\Psi(\Lambda)$: cf. lemmas \ref{lem-unsplit} and \ref{lem-extpsi} and their dual version.

\item The adjunction morphisms give natural epimorphisms 
$$i^{(k)}_* j^{(k)}_! \Lambda^{(k)}(\frac{1-k}{2}) \twoheadrightarrow 
\gr^k_!(\Psi(\Lambda))$$ 
where, when $\Lambda=\overline \Zm_l$, we never need to saturate
in the sense of definition \ref{defi-saturation}: cf. the discussion after the
proof of lemma \ref{lem-extpsi}.

\item Dually we also have monomorphisms 
$$\cogr_{*,k}(\Psi(\Lambda)) \hookrightarrow
i^{(r-k+1)}_* j^{(r-k+1)}_* \Lambda^{(r-k+1)}(\frac{r-k}{2})$$ 
induced by adjunction morphisms, where 
for $\Lambda=\overline \Zm_l$, the cokernel are torsion free.

\item For every $1 \leq h \leq r$ and for every $h \leq k \leq r$, the adjunction morphism
gives a epimorphism, cf. \ref{eq-psiFill1}, without taking saturation when $\Lambda=\overline \Zm_l$
$$i^{(k)}_* j^{(k)}_! \Lambda^{(k)}(\frac{k-1-2(h-1)}{2}) \twoheadrightarrow 
\Fill^{-k}(\gr^h_!(\Psi(\Lambda)).$$

\item Dually the adjunction morphism gives monomorphisms with torsion free cokernels, 
cf. \ref{eq-psiFill2}
$$\cogr_{*,h}(\Psi(\Lambda))/\Fill^k(\cogr_{*,h}(\Psi(\Lambda)) \hookrightarrow
i^{(k+1)}_* j^{(k+1)}_* \Lambda^{(k+1)}(\frac{1-k+2(h-1)}{2})$$
for $1 \leq h \leq r$ and $r-h \leq k \leq r-1$.
\end{itemize}

Combining compatibly the symbols $i_I^*$, $i_{I,*}$, $j_{I,*}$, $j_I^*$, the duality symbol $D$
and the versions $(h)$, on can compute their $\lexp p h^\delta$ either in the Grothendieck group
of perverse sheaves on $Y$ but also their filtrations and decide if the various extensions are split
or not. One example is given in the proposition \ref{prop-datlem} where we give a new proof of
\cite{dat-lem} theorem 2.2 in the strict semi stable case, cf. proposition \ref{prop-datlem}: 
in loc. cit. the author explains how to reduce
from the general semi stable case to the strict one.

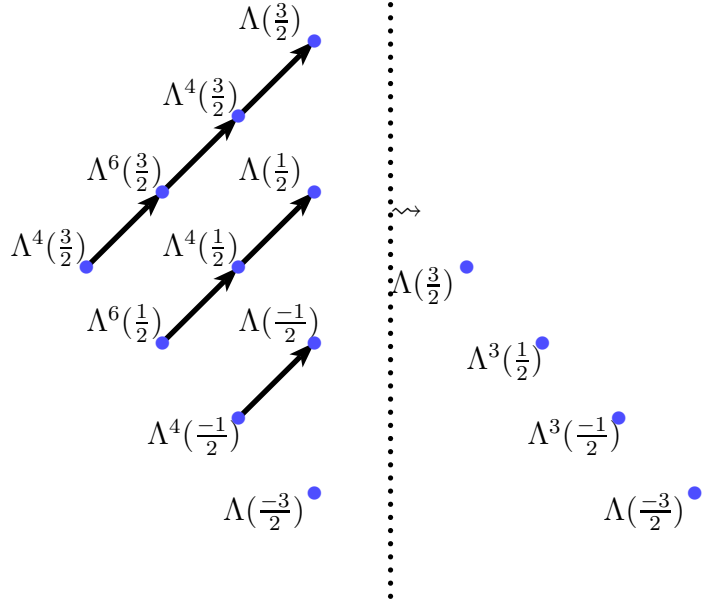
\begin{figure}[ht]
\centering
\input{fig-hi-rPsi.tex}
\caption{\label{fig-ss-hiPsi} Spectral sequence computing the stalks of $h^i \Psi$
at a geometric point of $Y_I^0$ with $\sharp I=4$: the abutment is on the right}
\end{figure}

The fact that the various extensions
inside $j_{I,!} \Lambda_I$ are non split, allows us, using the adjunction maps, to compute the sheaves
cohomology groups of everything, cf. figure (\ref{fig-ss-hiPsi}) and the diagram of the proposition
\ref {prop-ss-hiPsi}. It happens that the maps of the spectral sequence computing the
sheaf cohomology groups $h^i \Psi$, associated to the filtration constructed above, are in some
sense diagonal.

The main interest of these descriptions concerns the computation of the cohomology groups.
Usually, using trace formula, it is possible to obtain informations about the alternated sum
of the cohomology groups with compact supports. If one then wants to compute each individual
such cohomology groups, the strategy then consists, using the equality (\ref{eq-j!*}) to compute the alternated sum of the cohomology groups
of the intermediate extensions: purity then allows us to separate the individual cohomology groups.
It is moreover crystal clear to read the monodromy action. In \S \ref{para-application} we give
an illustration of such computations under the assumption that the strata behave as if they were
affine.

In a much more ramified situation given by the special fiber of Kottwitz-Harris-Taylor Shimura varieties
at a split prime $p$, in \cite{boyer-duke} we are also able to obtain the same results which at the
end looks very similar to the semi stable reduction case but the arguments are much more
complicated because first, contrary to the semi stable case, we do not know a priori the sheaf cohomology groups $R^i \Psi(\Lambda)$ and secondly the local systems which appears are not
constant sheaves so that their intermediate extensions for $p$ and $p+$ do not generally coincide.
To see how we exploit these descriptions to compute cohomology groups, cf. \cite{boyer-compositio}
or \cite{boyer-imj}.

\section{Nearby cycles in the semi-stable case}

Let denote by $K$ a finite extension of $\Qm_p$, with ring of integers $\OC_K$ and residue field
$\kappa$. Recall that a scheme $X \longrightarrow \spec \OC_K$ of finite type is said strictly
semi-stable purely of relative dimension $d-1$ if it is, Zariski locally on $X$, etale over 
$\spec \OC_K[T_1,\cdots,T_d]/(T_1 \cdots T_r - \varpi_K)$ for a prime element $\varpi_K$ in 
$\OC_K$ and an integer $1 \leq r \leq d$. Note that in particular  
\begin{itemize}
\item $X$ is regular and flat over $\spec \OC_K$, 

\item its generic fiber $X_\eta$ is smooth of relative dimension $d-1$,

\item its special fiber $Y:=X_s$ is a divisor of $X$ with simple normal crossings.
\end{itemize}
We now moreover suppose that $Y$ is globally the sum of smooth divisors $Y_i$ for
$i=1,\cdots,r$.

If $\overline K$ is a algebraic closure of $K$, and $X_{\overline \eta}:=X_\eta \times_{\spec K}
\spec \overline K$ with geometric special fiber $X_{\overline s}$, we have a diagram of
cartesian squares
$$\xymatrix{
X_{\overline \eta} \ar@{^{(}->}[r]^{\overline j} \ar[d] & \overline X \ar[d] & 
X_{\overline s} \ar@{_{(}->}[l]_{\overline i}  \ar[d] \\
\overline \eta \ar@{^{(}->}[r] & \overline S & \overline s \ar@{_{(}->}[l]
}$$
where $(\overline S,\overline \eta,\overline s)$ is the strict henselian triple attached to
$(S,\eta,s)$ and the choice of $\overline K$.

Let $\FC$ be a $\overline \Zm_l$-locally constant free sheaf on $X_\eta$ where $l \neq p$.
The complexe of nearby cycles of $\FC$
$$R\Psi(\FC):=\overline i^* R\overline j_* \FC_{\overline \eta}$$
is an object of the triangulated category of $\overline \Zm_l$-sheaves on $X_{\overline s}$ with
a continuous action of $\gal(\overline \eta/\eta)$ compatible with its action on $X_{\overline s}$:
moreover this action factors through the tame quotient $\gal(\eta_t/\eta)$.

Let $Y_1,\cdots,Y_r$ be the irreducible components of $Y:=X_{\overline s}$. For a non-empty subset
$I \subset \{ 1,\cdots,r \}$, we put $Y_I:=\bigcap_{i \in I} Y_i$ and let $i_I: Y_I \hookrightarrow Y$
be the associated immersion. The scheme $Y_I$ is smooth of dimension $d-k$ if $k=\sharp I$.
We then put
$$Y^{(k)}:=\bigcup_{\myatop{I \subset \{ 1,\cdots ,r \}}{\sharp I=k}} Y_I,$$
and let $i^{(k)}:Y^{(k)} \hookrightarrow Y$ be the associated immersion.
For every $I \subsetneq \{ 1,\cdots,r\}$, we put
$$Y_I^0:=Y_I \setminus \bigcup_{i \not \in I} Y_{I \cup \{ i \} },$$
and $j_I:Y_I^0 \hookrightarrow Y_I$ the open immersion. We also denote by
$$Y^{(k),0}:= \coprod_{\myatop{I \subset \{ 1,\cdots ,r \}}{\sharp I=k}} Y_I^0,$$
and $j^{(k)}:Y^{(k),0} \hookrightarrow Y^{(k)}$.

\begin{thm} (cf. \cite{ill} theorem 3.2 (c)) \\
Let $\Lambda$ be the constant sheaf
$\overline \Zm_l$ (or $\overline \Qm_l$ or $\overline \Fm_l$). We then have
\begin{itemize}
\item $R^0 \Psi(\Lambda) \simeq \Lambda$,

\item $R^1 \Psi(\Lambda) \simeq \bigl ( \bigoplus_{i=1}^r \Lambda_{Y_i} \bigr )/\Lambda ~(-1)$
where $\Lambda$ injects diagonally in the sum,

\item $R^q \Psi(\Lambda) \simeq \bigwedge^q R^1 \Psi(\Lambda)$.
\end{itemize}
\end{thm}

\section{Image of $\Psi(\Lambda)$ in the Grothendieck group}

We now denote by $\Psi(\Lambda):=R\Psi(\Lambda)[d-1](\frac{d-1}{2})$ which is then a perverse
sheaf of weight $0$. If $D$ denote the Grothendieck-Verdier duality functor then $D\Psi(\Lambda) \simeq
\Psi(\Lambda)$. In particular when $\Lambda=\overline \Zm_l$, then $\Psi(\Lambda)$ is a
torsion free perverse sheaf. 

The basic tool to determine the perverse irreducible constituants of some perverse sheaf whose
you only know its sheaf cohomology groups is given by the following classical statements.

\begin{prop}
Let $D_c^b(X,\Km)$ be the derived category of $\overline \Qm_l$-constructible complexes on a
$\Fm_q$-scheme $X$ equipped with the autodual perversity $t$-structure (resp. the
trivial $t$-structure) with core the category $\perv(X)$ (resp. $\const(X)$) of perverse sheaves
(resp. constructible sheaves) on $X$. For any object $\FC$ of $D_c^b(X,\overline \Qm_l)$, its
image in $\groth(\perv(X))$ is determined by its image in $\groth(\const(X))$.
\end{prop}

\begin{proof}
Given a locally small triangulated category $A$, we consider its Grothendieck group $K(A)$ defined
as the free group generated by the isomorphism classes of objects in $A$ quotiented by the
relations: $A=B+C$ for every distinguished triangle 
$B \longrightarrow A \longrightarrow C \xrightarrow{+1}$.
The result then follows from the following classical lemma.
\end{proof}

\begin{lem} 
Let $(\DC^{\leq 0},\DC^{\geq 0})$ be a derived category equipped with a
non-degenerate $t$-structure: we denote by $\CC$ its core, which
is then an abelian category of Grothendieck group $\groth(\CC)$. The map that associates to an object 
$\FC$ of $\DC$
$$\sum_i (-1)^i [\lexp p h^i \FC ] \in \groth(\CC)$$
induces an isomorphism of the Grothendieck group $K(\DC)$ of the triangulated category $\DC$ onto $\groth(\CC)$.
\end{lem}

\begin{proof} 
For any object $\GC$ of $\DC$ and for any $n \in \Zm$, we have a distinguished triangle
$$\tau_{\leq n} \GC \longrightarrow \GC \longrightarrow \tau_{\geq n+1} \GC \longmapright{+1}$$
such that if $\GC$ is an object of $\DC^{\geq n}$, we have $[\GC]=(-1)^n [\lexp p h^n \GC]+[\tau_{n+1} \GC]$
because $\tau_{\leq n} \GC=\tau_{\leq n} \tau_{\geq n} \GC=(\lexp p h^n \GC)[n]$. Let $a$ and $b$ be such that
$\FC$ is an object of $\DC^{[a,b]}$. Applying the above to $\tau_{\geq n} \FC$ for $n$ varying from
$a$ to $b$, we obtain the equality $[\FC]=\sum_i (-1)^i [\lexp p h^i \FC],$ so that the mapping $\sum_i
\alpha_i [P_i] \in \groth(\CC) \mapsto \sum_i \alpha_i [P_i] \in K(\DC)$ is the inverse of that in the statement.
\end{proof}

\begin{prop} \label{prop-perv1}
For $\Lambda=\overline \Qm_l$ or $\overline \Fm_l$, and
$Y$ a scheme equipped with a stratification
$$Y^{(r)} \subset Y^{(r-1 \leq r)} \subset \cdots \subset Y^{(1 \leq r)}=Y,$$
let $\PC$ be a perverse $\Lambda$-sheaf on $Y$. Then the image of $\PC$ in the Grothendieck group
of perverse $\Lambda$-sheaves on $Y$ is determined by those of $\sum_i (-1)^i [h^i(\PC)_{|Y^{(h)}}]$
in the Grothendieck group of locally constant sheaves on $Y^{(h)}:=Y^{((h \geq h+1)}-X^{(h+1)}$
for all $1 \leq h \leq r$.
\end{prop}

\begin{proof}
It suffices to consider the case of an open immersion $j:U \hookrightarrow Y$ of the complement
$i:Z \hookrightarrow Y$.
From the distinguished triangle 
$$j_! j^* \FC \longrightarrow \FC \longrightarrow i_*i^* \FC \xrightarrow{+1},$$ 
we deduce that
the image of an object $\FC$ of $D_c^b(Y,\Lambda)$ in $\groth(\const(Y))$ is determined by
that of $j^* \FC$ in $\groth(\const(U))$ and that of $i^* \FC$ in $\groth(\const(Z))$.
\end{proof}

%
We then want to describe the irreducible constituants of $\Psi(\Lambda)$
and its various filtrations by perverse sheaves.

\begin{nota} \label{nota-lambdaI}
For every $I \subset \{ 1,\cdots,r \}$, we put
$$\Lambda_I:=\Lambda_{Y_I^0}[d-\sharp I](\frac{d-\sharp I}{2}).$$
For $1 \leq h \leq d$, we put
$$\Lambda^{(h)}:=\Lambda_{Y^{(h),0}}[d-h](\frac{d-h}{2}).$$ 
All these shifted local systems are then perverse sheaves of weight $0$ on their support.
\end{nota}

\begin{lem}
For $\Lambda=\overline \Qm_l$ or $\overline \Fm_l$, and
for $I \subsetneq \{ 1,\cdots,r\}$, we have the following equality in the Grothendieck group of
perverse sheaves on $Y$:
\begin{equation} \label{eq-jI!}
i_{I,*} j_{I,!} \Lambda_I= \sum_{k=0}^{r-\sharp I} 
\sum_{\myatop{J \supset I}{\sharp J \setminus I=k}}
i_{J,*} \lexp p j_{J,!*} \Lambda_J(\frac{\sharp J - \sharp I}{2}).
\end{equation}
\end{lem}

\rem note that as the $j_{J}$ are affine, then $j_{I,!}$ preserve the perverse
$t$-structure. When $\Lambda=\overline \Zm_l$, it preserves both $t$-structures $p$ and $p+$
so that $j_{I,!} \Lambda_I$ are free perverse sheaves. Note also that, as 
$Y_I$ is smooth then $\lexp p j_{I,!*} \Lambda_I$ is equal to
$\Lambda_{Y_{I}}[d-\sharp I]$ and, when $\Lambda=\overline \Zm_l$, also equal to
$\lexp {p+} j_{I,!*} \Lambda_I$.

\begin{proof}
We start from $\lexp p j_{I,!*} \Lambda_I=\Lambda_{Y_I}[d-\sharp I]$ which by proposition 
\ref{prop-perv1} gives us the equality in the Grothendieck group
\begin{equation} \label{eq-j!*}
i_{I,*} \lexp p j_{I,!*} \Lambda_I=\sum_{k=0}^{r-\sharp I} 
\sum_{\myatop{J \supset I}{\sharp J \setminus I=k}} (-1)^{k} i_{J,*} 
j_{J,!} \Lambda_J(\frac{\sharp J-\sharp I}{2}).
\end{equation}
and so
\begin{equation} \label{eq-j!}
i_{I,*} j_{I,!} \Lambda_I=i_{I,*} \lexp p j_{I,!*} \Lambda_I + \sum_{k=1}^{r-\sharp I}
\sum_{\myatop{J \supset I}{\sharp J \setminus I=k}} (-1)^{k-1} i_{J,*}  j_{J,!} 
\Lambda_J(\frac{\sharp J-\sharp I}{2}).
\end{equation}
We then argue by induction on $\sharp I$ by assuming the result true for every $J$ such that
$\sharp J > \sharp I$: note that the case $\sharp J=r$ is obvious as 
$Y_{\{ 1,\cdots,r\} }^0=Y_{\{ 1,\cdots,r \} }$.
We then count the multiplicity of $i_{J,*} \lexp p j_{J,!*} \Lambda_J(\frac{\sharp J-\sharp I}{2})$ 
in the right hand side of
(\ref{eq-j!}) with $J=I \coprod \{ i_1,\cdots,i_\delta \}$. The contribution for a fixed $k$ in (\ref{eq-j!}),
by the induction hypothesis, is $\binom{\delta}{k} (-1)^{k-1}$ and the result follows from the classical
equality
$$1=\sum_{k=1}^{\delta} (-1)^{k-1} \binom{\delta}{k}.$$
\end{proof}

Considering all $I \subset \{ 1,\cdots, r \}$ with the same order,
the previous lemma become
\begin{equation} \label{eq-jh!*}
i^{(h)}_* \lexp p j^{(h)}_{!*} \Lambda^{(h)}=\sum_{k=0}^{r-h} (-1)^k
\binom{h+k}{h} i^{(h+k)}_* j^{(h+k)}_! \Lambda^{(h+k)}(\frac{k}{2})
\end{equation}
\begin{equation} \label{eq-jh!}
i^{(h)}_* j_{!}^{(h)} \Lambda^{(h)}= \sum_{k=0}^{r-h} \binom{h+k}{h} 
i^{(h+k)}_* \lexp p j^{(h+k)}_{!*} \Lambda^{(h+k)}(\frac{k}{2}).
\end{equation}

\begin{lem} \label{lem-important}
Let $s \in \{ 1,\cdots,r\}$ be fixed and put for every $1 \leq h \subset \leq r$
$$j^{(h)}_{\neq s}:Y^{(h)} \setminus (Y^{(h)} \cap Y_s ) \hookrightarrow Y^{(h)} 
\hookleftarrow Y^{(h)} \cap Y_s: i^{(h)}_{s}.$$
For $P:=\bigoplus_{\myatop{\sharp I=h}{s \not \in I}} \lexp p j_{I,!*} \Lambda_I$,
we have the following short exact sequence
$$0 \rightarrow i^{(h+1)}_*\lexp p j^{(h+1)}_{s,!*} \Lambda^{(h+1)}_s(\frac{1}{2}) \longrightarrow
i^{(h)}_* j^{(h)}_{\neq s,!} j^{(h),*}_{\neq s} P \longrightarrow P \rightarrow 0,$$
where $\Lambda_s^{(h+1)}$ is the sheaf $\Lambda$ supported on $\bigcup_{\myatop{\sharp I=h+1}{s \in I}}Y_I^0$. Moreover the extension does not split.
\end{lem}

\begin{proof}
By construction of intermediate extensions, note first that $\lexp pj^{(h)}_{\neq s,!*} j^{(h),*}_{\neq s} P
\simeq P$ so that $\lexp p h^0 i_s^{(h)} P=0$. As moreover the 
$\lexp p h^\delta i_s^{(h)} P=0$ for every $\delta<0$ to compute the last one, i.e. when $\delta=-1$.
we can simply first look at $\sum_\delta (-1)^\delta \lexp p h^\delta i_s^{(h)} P$ in
the Grothendieck group of perverse sheaves on $Y$. We start with the following equality
$$P=\sum_{k=0}^{r-h} \sum_{\myatop{\sharp J=h+k}{s \not \in J}} (-1)^k \binom{h+k}{h}
i_{J,*} \lexp p j_{J,!} \Lambda_J(\frac{k}{2}) + 
\sum_{k=1}^{r-h} \sum_{\myatop{\sharp J=h+k}{s \in J}} (-1)^k \binom{h+k-1}{h-1}
i_{J,*} \lexp p j_{J,!} \Lambda_J(\frac{k}{2}),$$
so that
$$\begin{array}{ll}
- i^{(h)}_* \lexp p h^{-1} i_s^{(h),*} P=\sum_\delta (-1)^\delta i^{(h)}_*
\lexp p h^\delta i_s^{(h),*} P & =
\sum_{k=1}^{r-h} \sum_{\myatop{\sharp J=h+k}{s \in J}} (-1)^k \binom{h+k-1}{h-1}
i_{J,*} \lexp p j_{J,!} \Lambda_J(\frac{k}{2}),\\
& = - i^{(h+1)}_* \lexp p j^{(h+1)}_{s,!*} \Lambda^{(h+1)}_s(\frac{1}{2}).
\end{array}$$
Put $j^{(h+1)}_s:\bigcup_{\myatop{\sharp I=h+1}{s \in I}} Y_I^0 \hookrightarrow Y^{(h+1)}$ and note that
$$j^{(h+1),*}_s \bigl ( \lexp p h^{-1} i_s^{(h)} \bigr ) \simeq \Lambda_s^{(h+1)}(\frac{1}{2})$$
so that the image of the adjunction morphism 
$$j^{(h+1)}_{s,*} j^{(h+1),*}_s  \bigl ( \lexp p h^{-1} i_s^{(h)} \bigr ) \longrightarrow
\lexp p h^{-1} i_s^{(h)}$$
has $\lexp p j^{(h+1)}_{s,!*} \Lambda^{(h+1)}_s(\frac{1}{2})$ as a quotient: the previous
equality in the Grothendieck group then tells us that it is an isomorphism.
\end{proof}

\begin{prop} \label{prop-psi-formula}
The image of $\Psi(\Lambda)$ in the Grothendieck group of perverse sheaves on $X$ is equal to
$$\sum_{h=1}^r \sum_{k=0}^{h-1} i^{(h)}_* \lexp p j^{(h)}_{!*} \Lambda^{(h)} (\frac{h-1-2k}{2})=
\sum_{i=0}^{r-1} \sum_{k=i+1}^{r} i^{(k)}_* \lexp p j^{(k)}_{!*} \Lambda^{(k)} (\frac{k-1-2i}{2}).$$
\end{prop}

\begin{proof}
We fix $I \subset \{ 1,\cdots,r\}$ and, for a closed geometric point $z$ of 
$Y_I^0$, we look at the stalk $(h^i \Psi)_z$ of $h^i \Psi$ at $z$. The answer is given by the
previous description of the $R^i \Psi(\Lambda)$ recalled above:
\begin{itemize}
\item for $i=0$,we obtain $\Lambda[d-1](\frac{d-1}{2})$,

\item for $i=1$, we have $\Lambda^{\sharp I-1}[d-2](\frac{d-2}{2})(-\frac{1}{2})$,

\item and for $i \geq 2$, the isomorphism $R^i \Psi(\Lambda) \simeq \bigwedge^i R^1 \Psi (\Lambda)$, 
gives
$$\Lambda^{\binom{\sharp I-1}{i}}[d-1-i](\frac{d-1-i}{2})(-\frac{i}{2}).$$
\end{itemize}

Now using (\ref{eq-jh!*}), for a fixed $0 \leq i \leq r-1$, 
the multiplicity of $i^{(h)}_* j^{(h)}_{!} \Lambda^{(h)} (\frac{h-1-2i}{2})$
in $\sum_{k=i+1}^{r} i^{(k)}_* \lexp p j^{(k)}_{!*} \Lambda^{(k)}(\frac{k-1-2i}{2})$ is equal to
$$\sum_{k=i+1}^{r} (-1)^{h-k} \binom{h}{k}=(-1)^{h-k-1} \binom{h-1}{i}.$$
Indeed starting with $\binom{h}{1}-\binom{h}{0}=\binom{h-1}{0}$, by an easy inductive argument we have
$$\binom{h}{k}-\binom{h}{k-1}+\cdots + (-1)^k \binom{h}{0}=
\binom{h}{k}- \binom{h-1}{k-2}=\binom{h-1}{k-1}.$$

We then deduce that for every $1 \leq h \leq r$, the restriction to $Y^{(h),0}$ of
$\sum_i (-1)^i h^i P$ for $P=\Psi(\Lambda)$ coincides with those for
$P$ equals to the formula given by the proposition. The result then follows from
the proposition \ref{prop-perv1}.
\end{proof}

\section{Filtrations of stratification from  \cite{boyer-torsion}: generalities}

Let $S$ be the spectrum of a finite field and $Y$ a finite-type scheme over $S$, then
the usual $t$-structure on $\DC(Y,\Lambda):=D^b_c(Y,\Lambda)$ is
$$\begin{array}{l}
A \in \lexp p \DC^{\leq 0}(Y,\Lambda)
\Leftrightarrow \forall x \in Y,~h^k i_x^* A=0,~\forall k >- \dim \overline{\{ x \} } \\
A \in \lexp p \DC^{\geq 0}(Y,\Lambda) \Leftrightarrow \forall x \in Y,~h^k i_x^! A=0,~\forall k <- \dim \overline{\{ x \} }
\end{array}$$
where $i_x:\spec \kappa(x) \hookrightarrow Y$ and $h^k(K)$ is the $k$-th sheaf of cohomology of $K$.

\begin{nota} 
We denote by $\lexp p \CC(Y,\Lambda)$ the heart pf these $t$-structure: the cohomological
functors are then denoted by $\lexp p h^i$. For a functor $T$, we put
$\lexp p T:=\lexp p hi^0 \circ T$.
\end{nota}

\rem $\lexp p \CC(Y,\Lambda)$ is a Noetherian abelian category and 
$\Lambda$-linear.
For $\Lambda$ a field, this $t$-structure is autodual for Verdier duality.
For $\Lambda=\overline \Zm_l$, we can equip the $\overline \Zm_l$-linear abelian category 
$\lexp p \CC(Y,\overline \Zm_l)$ with a torsion theory $(\TC,\FC)$ where 
$\TC$ (resp. $\FC$) is the full subcategory
of objects of $l^\oo$-torsion $T$ (resp. $l$-free $F$) , i.e. such that $l^N 1_T$ is zero for $N$ sufficiently large (resp. $l.1_F$ is a monomorphism).

\begin{defin} As in \cite{juteau}, we define the dual $t$-structure 
$$\begin{array}{l}
\lexp {p+} \DC^{\leq 0}(Y,\overline \Zm_l):= \{ A \in \lexp p \DC^{\leq 1}(Y,\overline \Zm_l):~
\lexp p h^1(A) \in \TC \} \\
\lexp {p+} \DC^{\geq 0}(Y,\overline \Zm_l):= \{ A \in \lexp p \DC^{\geq 0}(Y,\overline \Zm_l):~
\lexp p h^0(A) \in \FC \} \\
\end{array}$$
with heart $\lexp {p+} \CC(Y,\overline \Zm_l)$ equipped with its torsion theory
$(\FC,\TC[-1])$ dual of those of $\lexp p \CC(Y,\overline \Zm_l)$.
\end{defin}
%

\begin{nota} (cf. \cite{boyer-torsion} \S 1.3) \label{nota-FC}
Let
$$\FC(Y,\Lambda):=\lexp p \CC(Y,\Lambda) \cap \lexp {p+} \CC(Y,\Lambda)
=\lexp p \DC^{\leq 0}(Y,\Lambda) \cap \lexp {p+} \DC^{\geq 0}(Y,\Lambda)$$
be the quasi-abelian category of free perverse sheaves on $X$ with coefficients in $\Lambda$. 
\end{nota}

Let $j: U \hookrightarrow Y$ an affine open immersion and $i:F:=Y \setminus U \hookrightarrow Y$.
Then $j_!$ is $t$-exact and, cf. \cite{boyer-torsion} proposition 1.3.14, 
$j_!=\lexp p j_!=\lexp {p+} j_!$ with $j_! (\FC(U,\Lambda)) \subset \FC(Y,\Lambda)$. Starting from
the distinguished triangle
$$j_!j^* L \longrightarrow L \longrightarrow i_*i^* L \leadsto,$$
and using the perversity of $L$ and $j_! j^* L$, the long exact sequence associated to it gives
$$0 \rightarrow i_* \lexp p h^{-1} i^* L \longrightarrow \lexp p j_! j^* L
\longrightarrow L \longrightarrow i_* \lexp p h^0i^* L \rightarrow 0,$$
where $ i_* \lexp p h^{-1} i^* L$ is free as $\lexp p j_! j^* L=j_!j^* L$ is free.

\begin{defin} \label{defi-strict}
A monomorphism $f:L \hookrightarrow L'$ in $\lexp p \CC(Y,\Lambda)$ between free perverse
sheaves, is said strict if, and only if, its cokernel in $\lexp p \CC(Y,\Lambda)$ is free.
This is equivalent to ask $f$ to be a monomorphism of 
$\lexp {p+} \CC(Y,\Lambda)$.
\end{defin}

\begin{defin} A \emph{bimorphism} of $\FC(Y,\Lambda)$ is a monomorphism which is also
a epimorphism, and we will denote by $L \htarrow L'$ such a bimorphism between $L$ and $L'$. 
If moreover the kernel in $\lexp {p+} \CC(Y,\Lambda)$ has its support of dimension strictly less that
the dimension of the support of $L$, we will denote by $L \htarrow_+ L'$.
\end{defin}

For $L \in \FC(Y,\Lambda)$, we consider the following diagram 
$$\xymatrix{ 
& L \ar[drr]^{\can_{*,L}} \\
\lexp {p+} j_! j^* L \ar[ur]^{\can_{!,L}} 
\ar@{->>}[r]|-{+} & \lexp p j_{!*}j^* L \ar@{^{(}->>}[r]_+ & 
\lexp {p+} j_{!*} j^* L \ar@{^{(}->}[r]|-{+} & \lexp p j_*j^* L
}$$
where the bottom line is, cf. the remark following 1.3.12 in \cite{boyer-torsion},
is the canonical factorization of 
$\lexp {p+} j_! j^* L \longrightarrow \lexp {p} j_* j^* L$
and where the maps  $\can_{!,L}$ and $\can_{*,L}$ are given by adjunction.

\begin{nota} \label{nota-filtration0}
We put
$\PC_L:=i_*\lexp p h^{-1}_{libre} i^*j_* j^* L =\ker_\FC \Bigl ( \lexp {p+} j_! 
j^* L \twoheadrightarrow \lexp p j_{!*} j^* L \Bigr ).$
With the notations of the above diagram, we put
$$\Fil^0_{U,!}(L)=\im_\FC (\can_{!,L}) \quad \hbox{ et } \quad \Fil^{-1}_{U,!}(L)=
\im_\FC \Bigl ( (\can_{!,L})_{|\PC_L} \Bigr ).$$
\end{nota}

We now suppose that $Y$ is equipped with a stratification
$$Y^{(r)} \subset Y^{(r-1)} \subset \cdots \subset Y^{(1)}.$$
For $1 \leq h <r$, we put $Y^{(1 \leq h)}:=Y^{(1)}-Y^{(h+1)}$ and
$j^{(1 \leq h)}:Y^{(1 \leq h)} \hookrightarrow Y^{(1)}$.

\begin{defin} \label{defi-fil!}
For $P \in \FC(Y,\Lambda)$ we define 
$$\Fil^h_{!}(L):=\im_\FC \Bigl ( \lexp {p} j^{(1 \leq h)}_! j^{(1 \leq h),*} P \longrightarrow P \Bigr ).$$
\end{defin}

\begin{prop} \label{prop-ss-filtration} (cf. \cite{boyer-torsion} \S 2.2)
The previous definition functionally equips every $L \in \FC(Y,\Lambda)$ 
with a filtration of stratification.
$$0=\Fil^{0}_{!}(L) \subset \Fil^1_{!}(L) \subset \Fil^2_{!}(L) \cdots \subset 
\Fil^{r-1}_{!}(L) \subset \Fil^r_{!}(L)=L.$$
\end{prop}

\begin{defin} \label{defi-saturation}
For $L,L' \in \FC(Y,\overline \Zm_l)$ and $f:L \longrightarrow L'$ a morphism
the usual cokernel $L/\im_\CC f$
where $\im_\CC f$ is the image of $f$ in the category $\lexp p\CC(Y,\Lambda)$, might have non trivial
torsion. However if we consider $L/\im_{\CC^+} f$ where $\im_{\CC^+} f$ is the image of
$f$ in $\lexp {p+} \CC(Y,\Lambda)$, then it is necessary free. We call this modification
\emph{the saturation process}.
\end{defin}

\rem When $Y$ is a point and $f: L:=\Zm_l \longrightarrow L=\Zm_l$ is the multiplication by $l$,
then $L/\im_\CC f \simeq \Zm_l/l\Zm_l$ while $L/\im_{\CC^+} f \simeq \Zm_l/\Zm_l$. In other words
$\im_{\CC^+} f= \Zm_l$ is the usual saturation of $l \Zm_l \subset \Zm_l$.

The above filtrations of stratification are generally not fine enough. 
In \cite{boyer-torsion} proposition 2.3.3, using $\Fil^{-1}_{U,!}$, we construct
in a functorial manner the exhaustive filtration of stratification of any object $L$ of $\FC(Y,\Lambda)$, 
\begin{equation} \label{eq-fill}
0=\Fill^{-2^{r-1}}_{!}(L) \subset \Fill^{-2^{r-1}+1}_{!}(L) \subset \cdots \subset
\Fill^0_{!}(L) \subset \cdots \subset \Fill^{2^{r-1}-1}_{!}(L)=L,
\end{equation}
such that all graduates $\grr^k$ are free and are, on $\overline \Qm_l$,
simple, i.e. verify $\lexp p j^{(h)}_{!*} j^{(h),*} \grr^k  \htarrow_+ \grr^k$,
where $Y^{(h)}$ is the support of $\grr^k$.

Dually using $\can_{*,L}$ with
$$\coFil_{*,h}(L):=\coim_{\FC} \bigl ( L \longrightarrow j_*^{(h)} j^{(h),*} L \bigr ),$$
we also define a cofiltration
$$L=\coFil_{*,r}(L) \twoheadrightarrow \coFil_{*,r-1}(L) \twoheadrightarrow \cdots 
\twoheadrightarrow \coFil_{*,1}(L) \twoheadrightarrow \coFil_{*,0}(L)=0,$$
and a filtration
$$0=\Fil_{*}^{-r}(L) \hookrightarrow \Fil_*^{-r+1}(L) \hookrightarrow \cdots \hookrightarrow
\Fil_*^{-1}(L) \hookrightarrow \Fil_*^0(L)=L$$
where $\Fil_*^{-h}(L):=\ker_\FC \bigl ( L \twoheadrightarrow \coFil_{*,h}(L) \bigr )$.

\section{Filtrations of $j_{I,!} \Lambda_I$}
\label{para-fil!}

For $\Lambda=\overline \Qm_l$,
the weight filtration of $j_{I,!} \Lambda_I$, thanks to the equality (\ref{eq-jI!}),
can be written, up to a shift of $\sharp I$:
$$(0)=\Fil^W_{-r-1}(I,!) \subset \Fil^W_{-r}(I,!) \subset \cdots \subset \Fil^W_{\sharp I}(I,!)
=j_{I,!} \Lambda_I$$
with graded parts $\gr^W_{-k}(I,!)=\sum_{\myatop{J \supset I}{\sharp J=k}}
i_{J,*} \lexp p j_{J,!*} \Lambda_J(\frac{\sharp J-\sharp I}{2}).$
For $\Lambda=\overline \Zm_l$, we also have a similar filtration. Moreover as
the intermediate extensions for the two $t$-structures $p$ and $p+$ are the same for
the local systems $\Lambda_{Y_I^0}$ so that for $\Lambda=\overline \Fm_l$ we also have
a filtration as above with graduate parts the $\sum_{\myatop{J \subset I}{\sharp J \setminus I=k}}
i_{J,*} \lexp p j_{J,!*} \Lambda_J(\frac{\sharp J-\sharp I}{2}).$

\begin{lem} \label{lem-unsplit}
For $\Lambda=\overline \Qm_l$ or $\overline \Fm_l$ and
for every $\sharp I < k < r-1 $, the extension
$$0 \rightarrow \gr^W_{-k}(I,!) \longrightarrow \Fil^W_{-k+1}(I,!)/\Fil^W_{-k-1}(I,!) \longrightarrow
\gr^W_{-k+1}(I,!) \rightarrow 0,$$
does not split. 
\end{lem}

\begin{proof}
Note first that this is obvious for $k=\sharp I+1$ because 
$i_{I,*} \lexp p j_{I,!*} \Lambda_{Y_I^0}[d-\sharp I]$
is the top of $i_{I,*} j_{I,!} \Lambda_{Y_I^0}[d-\sharp I]$. We now suppose $k>\sharp I+1$ and 
we look at the spectral sequence 
$$E_1^{p,q}=h^{p+q} \gr^W_{-p}(I,!) \Rightarrow h^{p+q} j_{I,!} \Lambda_{Y_I^0}[d-\sharp I],$$
computing the sheaf cohomology groups of the extension by zero of $\Lambda_{Y_I^0}[d-\sharp I]$.
Consider any $J \supset I$ with $\sharp J=-k$ and a generic point $z$ of $Y_J^0$.
Concerning the stalks at $z$:
\begin{itemize}
\item $(E_\oo^n)_z$ is zero,

\item $(h^{p+q} \gr^W_{-k}(I,!))_z$ is non zero if and only $p+q=d-\sharp J$,

\item $(h^{p+q} \gr^W_{-k+1}(I,!))_z$ is non zero if and only $p+q=d-\sharp J-1$,

\item for every $\delta \not \in \{ -k,-k+1 \}$,
$(h^{p+q} \gr^W_{\delta}(I,!))_z=0$ for $p+q \in \{ d-\sharp J-1, d-\sharp J, d-\sharp J+1 \}$.
\end{itemize}
Thus for any filtration of $j_{I,!} \Lambda_{Y_I^0}[d-\sharp I]$, the index of the graded part
containing $\gr^W_{-k}(I,!))_z$ should be less than those containing
$\gr^W_{-k+1}(I,!))_z$, which proves the statement.
\end{proof}

\rem otherwise stated, up to obvious modifications such as refinement or forgetting some terms, 
this weigh filtration of $j_{I,!} \Lambda_I$ is unique. As a consequence we deduce that the
adjunction morphism
\begin{equation} \label{eq-jFill1}
\bigoplus_{\myatop{J \supset I}{\sharp J=k}} j_{J,!} \Lambda_J(\frac{\sharp J-\sharp I}{2})
\simeq j^{(k)}_! j^{(k),*} \Fil^W_{-k}(I,!) \longrightarrow \Fil^W_{-k}(I,!),
\end{equation}
is surjective for $\Lambda=\overline \Qm_l$ and $\overline \Fm_l$. In particular it is also
surjective for $\Lambda=\overline \Zm_l$ without saturating the image. Dually the weight
filtration of $i_{I,*} j_{I,*} \Lambda_I$ up to a shift of $\sharp I$ is
$$(0)=\Fil^W_{r+1}(I,*) \subset \Fil^W_{r}(I,*) \subset \cdots \subset \Fil^W_{\sharp I}(I,*)=
j_{I,*} \Lambda_I,$$
with graded parts $\gr^W_k(I,*)=\sum_{\myatop{J \supset I}{\sharp J\setminus I=r-k}}
i_{J,*} \lexp p j_{J,!*} \Lambda_J(\frac{\sharp I-\sharp J}{2})$.
Moreover the adjunction morphism
\begin{equation} \label{eq-jFill2}
j_{I,*}\Lambda_I / \Fill^W_{k+1}(I,*) \hookrightarrow 
j^{(r-k+\sharp I)}_! j^{(r-k+\sharp I),*} \Bigl (j_{I,*}\Lambda_I / \Fill^W_{k+1}(I,*) \Bigr )
\simeq \bigoplus_{\myatop{J \supset I}{\sharp J\setminus I=r-k}}
j_{J,*} \Lambda_J(\frac{\sharp I-\sharp J}{2}),
\end{equation}
is injective with, for $\Lambda=\overline \Zm_l$, a torsion free cokernel.

\begin{nota} Let 
$P_I=i_{I,*} \lexp p h^{-1} i_I^* j_{I,*} \Lambda_I)$ be the kernel of the canonical
morphism
$$P_I=\ker \Bigl (  j_{I,!} \Lambda_I \longrightarrow
\lexp p j_{I,!*} \Lambda_I \Bigr ).$$
We denote by $Y_{I^+}:=\bigcup_{J \supsetneq I} Y_J$ its support and as before
$j_{I^+}:Y_{I^+}^0 \hookrightarrow Y_{I^+}$.
\end{nota}

From the previous lemma, we deduce that, for 
$\Lambda \in \{ \overline \Qm_l,\overline \Zm_l,\overline \Fm_l \}$, the adjunction morphism
$$j_{I^+,!} j_{I^+}^* P_I \twoheadrightarrow P_I,$$
is surjective. Indeed the cokernel has support in $Y_I \setminus Y_I^0$ and the previous lemma
tells us that, over $\overline \Qm_l$ and $\overline \Fm_l$, 
the top of $P_I$ is $\lexp p j_{I,!*} j_I^* P_I$ so that the cokernel is zero. In particular,
for $\Lambda=\overline \Zm_l$, there is no need
to use a saturation process as in definition \ref{defi-saturation}.

\begin{lem} \label{lem-jneq!}
For $J_0 \supset I$ with $\sharp J_0 \setminus I=1$, we have a short exact sequence
$$0 \rightarrow j_{I^+ \setminus J_0,!} j_{I^+ \setminus J_0}^* P_I \longrightarrow P_I
\longrightarrow \lexp p j_{J_0,!*} \Lambda_J(\frac{1}{2}) \rightarrow 0,$$
and a epimorphism
$$j_{I^+,!} j_{I^+}^* P_I \twoheadrightarrow P_I$$
without taking saturation.
\end{lem}

\begin{proof}
By construction $P_I$ (resp. $j_{I^+ \setminus J_0,!} j_{I^+ \setminus J_0}^* P_I$) 
is the constant sheaf $\Lambda[d-\sharp I-1](\frac{d-\sharp I}{2})$ on $Y_{I^+}$ 
(resp. on $Y_{I^+} \setminus Y_{J_0}$)
while $\lexp p j_{J_0,!*} \Lambda_J$ is the constant sheaf
$\Lambda[d-\sharp J_0](\frac{d-\sharp J_0}{2})$ on $Y_{J_0}$ with $\sharp J_0=\sharp I +1$.

The second statement is proved in the same way.
\end{proof}

Dually we obtain we obtain a cofiltration of $L:=j_{I,*} \Lambda_I$:
$$i_{I,*}j_{I,*} \Lambda_I= \coFil_{*,r-\sharp I+1}(L) \twoheadrightarrow \coFil_{*,r-\sharp I}(L) 
\twoheadrightarrow \cdots \twoheadrightarrow \coFil_{*,1}(L) \twoheadrightarrow \coFil_{*,0}(L)=0$$
with graded parts $\cogr_{*,k}(L):=\ker(\coFil_{*,k}(L) \twoheadrightarrow \coFil_{*,k-1}(L) )$ 
having image in the Grothendieck group of perverse sheaves
on $X$ equal to
$$\sum_{J \supset I} i_{J,*} \lexp p j_{J,!*} \Lambda_J(\frac{\sharp I - \sharp J}{2}),$$
the extension between two consecutive graded parts being moreover unsplit.

\begin{prop} \label{prop-ij}
For every $I \subset J$ we have
$$i_{J,*} \lexp p h^0 i_J^* \bigl ( j_{I,*} \Lambda_I \bigr ) \simeq i_{I,*} j_{I,*} \Lambda_J(-\frac{1}{2}).$$
\end{prop}

\begin{proof}
By an obvious induction, we only need to deal with the case $\sharp J \setminus I=1$.
Put $j_{\overline I \setminus J}:Y_I \setminus Y_J \hookrightarrow Y_I$ and recall that
$i_{J,*} \lexp p h^0 i_J^* \bigl ( j_{I,*} \Lambda_I \bigr )$ is the cokernel of the adjunction morphism
$$j_{\overline I \setminus J,!} j_{\overline I \setminus J}^* j_{I,*} \Lambda_I \longrightarrow
j_{I,*} \Lambda_I,$$
so that, considering the constituant of $j_{I,*} \Lambda_I$ having support in $Y_J$, 
its image in the Grothendieck group of perverse sheaves on $Y_I$ is less or equal to
$$\sum_{H \supset J} \lexp p j_{H,!*} \Lambda_H (\frac{\sharp I - \sharp H}{2}).$$
This is in fact an equality as, (\ref{eq-jI!}) tells us that any of the 
$\lexp p j_{H,!*} \Lambda_H (\frac{\sharp I - \sharp H}{2})$ is not a subquotient of 
$j_{\overline I \setminus J,!} j_{\overline I \setminus J}^* j_{I,*} \Lambda_I$.

In particular the cofiltration of
stratification of $L:=i_{J,*} \lexp p h^0 i_J^* \bigl ( j_{I,*} \Lambda_I \bigr )$
$$ \coFil_{*,r-\sharp J+1}(L) \twoheadrightarrow \coFil_{*,r-\sharp J}(L) 
\twoheadrightarrow \cdots \twoheadrightarrow \coFil_{*,1}(L) \twoheadrightarrow \coFil_{*,0}(L)=0$$
is such that its graded parts $\cogr_{*,k}(L)$ have image in the Grothendieck group of perverse sheaves
on $X$ equal to
$$\sum_{H \supset J} \lexp p j_{H,!*} \Lambda_H (\frac{\sharp I - \sharp H}{2}),$$
where moreover the extension between two consecutive graded parts are unsplit.
The same being true for the cofiltration of stratification of $j_{J,*} \Lambda_J(-\frac{1}{2})$
we then deduce that the adjunction morphism
$$L \longrightarrow j_{J,*} j_J^* L$$
is a isomorphism.
%
\end{proof}

\section{Filtrations of $\Psi(\Lambda)$}

Let 
$$\Fil_!^p(\Psi(\Lambda))=\im_\FC \Bigl (j^{(1 \leq p)}_! j^{(1 \leq p),*} \Psi(\Lambda) 
\longrightarrow \Psi(\Lambda) \Bigr ),$$
be the filtration of stratification of $\Psi(\Lambda)$ 
$$(0)=\Fil_!^0(\Psi) \subset \Fil_!^1(\Psi) \subset \cdots \subset \Fil_!^r(\Psi)=\Psi(\Lambda).$$

\begin{prop}
For $1 \leq k \leq r$, the image of $\gr_!^k(\Psi(\Lambda)) \otimes_\Lambda \overline \Qm_l$ 
in the Grothendieck group 
of perverse  $\overline \Qm_l$-sheaves on $X$ is equal to that of
$\sum_{i=k}^{r} i^{(i)}_* \lexp p j^{(i)}_{!*} \Lambda^{(i)} (\frac{i-1-2(k-1)}{2}).$
Moreover $\gr_!^k(\Psi(\Lambda))$ has a filtration 
$$(0)=\Fill^{-r-1}(\gr_!^k(\Psi(\Lambda))) \subset \Fill^{-r}(\gr_!^k(\Psi(\Lambda))) \subset \cdots \subset
\Fill^{-k}(\gr_!^k(\Psi(\Lambda)))=\gr_!^k(\Psi(\Lambda))$$
with $\grr^{-h}(\gr_!^k(\Psi(\Lambda))) \simeq  i^{(h)}_*
\lexp p j^{(h)}_{!*} \Lambda^{(h)} (\frac{h-1-2(k-1)}{2})$.
\end{prop}

\rem 
Note first that, considering the weight in the equality of proposition \ref{prop-psi-formula}, the image of 
$\gr_!^k(\Psi(\Lambda)) \otimes_\Lambda \overline \Qm_l$ is necessary less than
$\sum_{i=k}^{r} i^{(i)}_* \lexp p j^{(i)}_{!*} \Lambda^{(k)} (\frac{i-1-2(k-1)}{2}).$
Moreover as explained in the previous section, lemma \ref{lem-unsplit} implies that the filtration
$\Fill_\bullet$ of each $\gr_!^k(\Psi(\Lambda))$ is essentially unique.

\begin{proof}
We argue as in the proof of lemma \ref{lem-unsplit} by proving the following observation.
Consider $I$ with $\sharp I=i+1$ then there does not exists a filtration of
$\Psi(\Lambda)$ such that there exists two graduate parts respectively isomorphic to 
$\lexp p j_{I,!*} \Lambda^{(k)} (\frac{i-2(k-1)}{2})$ and
$\lexp p j^{(i)}_{!*} \Lambda^{(i)} (\frac{i-1-2(k-1)}{2})$ with respective index $r_1 < r_2$.

Indeed otherwise for $z$ a generic point of $Y_I^0$ the stalks at $z$ with weight 
$2(k-1)-i+1$ appears 
\begin{itemize}
\item in degree $\sharp I-d$ for $\lexp p j_{I,!*} \Lambda^{(k)} (\frac{i-2(k-1)}{2})$,

\item in degree $\sharp I-d-1$ for $\lexp p j^{(i)}_{!*} \Lambda^{(i)} (\frac{i-1-2(k-1)}{2})$

\item in degree $< \sharp I-d-1$ otherwise.
\end{itemize}
So if $r_1<r_2$, we then see that $h^{-d+\sharp I}\Psi(\Lambda)$ would have a non zero stalks at $z$
with weight $2(k-1)-i+1$ which is not the case.
We then deduce that every irreducible constituant of 
$\sum_{i=k}^{r} i^{(i)}_* \lexp p j^{(i)}_{!*} \Lambda^{(i)} (\frac{i-1-2(k-1)}{2})$ has to be a constituant
of $\Fil_!^k(\Psi)$ which gives the result.
\end{proof}

As the $p$ and $p+$ intermediate extensions of our local systems are isomorphic, we then
obtain for $\Lambda=\overline \Fm_l$, a filtration like in the previous proposition
with graded parts the $\lexp p j^{(h)}_{!*} \Lambda^{(h)} (\frac{h-1-2(k-1)}{2})$.
Consider then two consecutive graded parts in the filtration $\Fill^\bullet (\gr^k_!(\overline \Fm_l))$,
say for example $P_I:=\lexp p j_{I,!*} \Lambda_I(\frac{\sharp I-1-2(k-1)}{2})$ and
$P_-:=\lexp p j_{J,!*} \Lambda_J(\frac{\sharp J-1-2(k-1)}{2})$ with $\sharp I=\sharp J+1$. We choose
$s \in I \setminus J$ and consider the extension $X$ inside $\gr^k_!(\overline \Fm_l)$:
$$0 \rightarrow P_I \longrightarrow X \longrightarrow P_- \rightarrow 0.$$

\begin{lem} \label{lem-extpsi}
The above extension $X$ does not split.
\end{lem}

\begin{proof}
We first start by proving the following statement which
the analog of the lemma \ref{lem-jneq!} for $\Psi(\Lambda)$.

\begin{lem} \label{lem-analogpsi}
For $h \in \{ 1,\cdots,r \}$, we put $j_{\neq h}:Y \setminus Y_h \hookrightarrow Y$. Then the cokernel
$L:=i_{h,*} \lexp p h^0 i_h^* \Psi(\Lambda)$ of the adjunction morphism
$$j_{\neq h,!} j_{\neq h}^* \Psi(\Lambda) \longrightarrow \Psi(\Lambda)$$
is isomorphic to $j_{h,*} \Lambda^{(h)}$.
\end{lem}

\begin{proof}
Recall the short exact sequence
$$0 \rightarrow \Psi(\Lambda) \longrightarrow \overline j_! \Lambda[d-1](\frac{d-1}{2}) 
\longrightarrow \overline j_* \Lambda[d-1](\frac{d-1}{2}) \rightarrow 0,$$
with $\Psi(\Lambda)=\lexp p h^{-1} \overline i^* \overline j_* \Lambda[d-1](\frac{d-1}{2})$.
Note also that for every $\delta \neq -1$, we have
$\lexp p h^{\delta} \overline i^* \overline j_* \Lambda[d-1]=(0)$.
As $\overline j_{\neq h}: X_{\bar \eta} \setminus Y_h \hookrightarrow X_{\bar \eta}$ is 
affine,
for $\overline i_h: Y_h \hookrightarrow X_{\bar \eta}$, the long exact sequence associated
to the distinguished triangle
$$\overline j_{\neq h,!} \overline j_{\neq h}^* \Lambda[d-1] \longrightarrow \Lambda[d-1]
\longrightarrow \overline i_{h,*} \overline i_h^* \Lambda[d-1] \leadsto^{+1},$$
as the first two terms are perverse, then
$\lexp p h^{\delta} \overline i_h^* \Lambda[d-1]=(0)$ for $\delta \neq-1$. Writting
$\overline i_h=\overline i \circ i_h$ and considering the spectral sequence
$$E_2^{r,s}=\lexp p h^r i_h^* \lexp p h^s \overline i^* \Lambda[d-1]
\Rightarrow \lexp p h^{r+s} \overline i_h^* \Lambda[d-1],$$
we then deduce that $i_h^* \Psi(\Lambda)$ is perverse and that
$$0 \rightarrow j_{\neq h,!} j_{\neq h}^* \Psi(\Lambda) \longrightarrow \Psi(\Lambda) \longrightarrow
i_{h,*} \lexp p h^0 i_h^* \Psi(\Lambda) \rightarrow 0.$$

We then observe that the stalks at geometric points of 
$\Psi(\Lambda)$ and that of 
$$\bigoplus_{\myatop{J \supset \{ h \} }{\sharp J=r-k+1}} \lexp p j_{J,!*} \Lambda_J (-\frac{r-k}{2}),$$
are the same, so that the image of $i_{h,*} \lexp p h^0 i_h^* \Psi(\Lambda)$ in the Grothendieck
group of perverse sheaves on $Y$ is equal to that of this direct sum. Moreover
$\Psi(\Lambda)$ has a quotient $Q$ with the following cofiltration of stratification
$$Q=\coFil_{*,r}(Q) \twoheadrightarrow \coFil_{*,r-1}(Q) \twoheadrightarrow \cdots \twoheadrightarrow
\coFil_{*,1}(Q) \twoheadrightarrow \coFil_{*,0}(Q)=0,$$
with graded parts $\cogr_{k,*}(Q)$ having image in the Grothendieck group of perverse sheaves on $X$
equal to
$$\sum_{\myatop{J \supset \{ h \} }{\sharp J=r-k+1}} i_{J,*} \lexp p j_{J,!*} \Lambda_J (-\frac{r-k}{2}),$$
where the extensions are unsplit. We then deduce that
$Q$ is isomorphic with $L$ and that the adjunction morphism
$$L \longrightarrow j_{h,*}j_h^* L$$
is an isomorphism.
\end{proof}

We now go back to the extension $X$ and we consider the short exact sequence
$$0 \rightarrow j_{\neq s,!} j_{\neq s}^* \Psi(\Lambda) \longrightarrow \Psi(\Lambda) \longrightarrow
i_{h,*} \lexp p h^0 i_s^* \Psi(\Lambda) \rightarrow 0,$$
so that our extension $X$ lives inside the first term $j_{\neq s,!} j_{\neq s}^* \Psi(\Lambda)$
By $t$-exactness of $j_{\neq s,!} j_{\neq s}^*$ and using lemma \ref{lem-important} we see
that $j_{\neq s,!} j_{\neq s}^* \Psi(\Lambda)$ can be written, up to a twist of weights, 
as successive unsplit extensions of the form
$$0 \rightarrow \lexp p j^{(h+1)}_{s,!*} \Lambda^{(h+1)}_s(\frac{1}{2}) \longrightarrow
j^{(h)}_{\neq s,!} j^{(h),*}_{\neq s} \bigl ( \bigoplus_{\myatop{\sharp I=h}{s \not \in I}} \lexp p j_{I,!*}
\Lambda_I \bigr ) \longrightarrow \bigoplus_{\myatop{\sharp I=h}{s \not \in I}} \lexp p j_{I,!*}
\Lambda_I  \rightarrow 0,$$
in which $X$ appears.
\end{proof}

In other words the adjunction morphism giving rise to
$$j^{(k)}_! \Lambda^{(k)}(\frac{1-k}{2}) \twoheadrightarrow \gr^k_!(\Psi(\Lambda))$$
is surjective when $\Lambda=\overline \Zm_l$, i.e. without using the saturation process of
definition \ref{defi-saturation}. 
As a consequence, using (\ref{eq-jFill1}), the adjunction morphisms give epimorphism
\begin{equation} \label{eq-psiFill1}
j^{(k)}_! \Lambda^{(k)}(\frac{k-1-2(h-1)}{2}) \twoheadrightarrow 
\Fill^{-k}(\gr^h_!(\Psi(\Lambda))),
\end{equation}
whose kernel $K(k,h)$ has a filtration
$$\Fill^{-r-1}(K(k,h))=(0) \subset \Fill^{-r}(K(h,k)) \subset \cdots \subset \Fill^{-k-1}(K(k,h))=K(k,h)$$
with graded parts $\grr^{-i} K(h,k) \simeq (\binom{i}{k}-1) \lexp p j^{(i)}_{!*} \Lambda^{(i)}
(\frac{i-1-2(h-1)}{2})$.  Note also that the image of (\ref{eq-psiFill1}) in the Grothendieck group
of perverse sheaves on $Y$ is equal to $\gr^k_!(\Psi(\Lambda))(\frac{k-h}{2})$ whose is also a quotient
of $j^{(k)}_! \Lambda^{(k)}(\frac{k-1-2(h-1)}{2})$ so that (\ref{eq-psiFill1}) induces a strict monomorphism
$$\gr^k_!(\Psi(\Lambda))(\frac{k-h}{2}) \hookrightarrow \gr^h_!(\Psi(\Lambda)).$$
Moreover all the process to construct filtrations never needs
saturation when $\Lambda=\overline \Zm_l$.

\begin{prop} \label{prop-ss-hiPsi}
Let $I \subset \{ 1,\cdots,r\}$ non empty and consider a geometric point $z_I$ of
$Y_I^0$. The spectral sequence computing the stalks at $z_I$ of the $h^i \Psi(\Lambda)$
looks like as in figure \ref{fig-ss-hiPsi} and degenerates at $E_2$. More precisely all the maps
$h^i \gr_!^p(\Psi(\Lambda) \longrightarrow h^{i+1} \gr_!^{p-\delta}(\Psi(\Lambda))$ for any
$\delta >0$ are zero while the maps
$$h^i \grr^{-p}(\gr^k_!(\Psi(\Lambda))) \longrightarrow h^{i+1} \grr^{-p-\delta}(\gr^k_!(\Psi(\Lambda)))$$
are non zero only for $\delta=1$ and given by a linear map
$$\bigoplus_{\myatop{J \subset I}{\sharp J=p}} \lexp J \Lambda \longrightarrow 
\bigoplus_{\myatop{K \subset I}{\sharp K=p+1}} \lexp K \Lambda$$
where the notation $\lexp J \Lambda$ and $\lexp K \Lambda$ represents a copy of 
$\Lambda(\frac{\sharp I-2k+1}{2})$,
and the matrix of this linear map in the canonical basis is $(\delta_{J,K})$ where
$\delta_{J,K}=1$ is $J \subset K$ and zero otherwise
$$\xymatrix{
z_I^* h^{k-r} \Psi(\Lambda) \ar@{^{(}->}[r]  \ar[d]^\simeq &
\bigoplus_{\myatop{J_k \subset I}{\sharp J=k}} \lexp {J_k} \Lambda \ar[d]^\simeq \ar[r] &
\bigoplus_{\myatop{J_{k+1} \subset I}{\sharp J=k+1}} \lexp {J_{k+1}} \Lambda \ar[d]^\simeq \ar[r] &
\cdots \ar[r] & 
\bigoplus_{\myatop{J_{\sharp I-1} \subset I}{\sharp J=\sharp I-1}} \lexp {J_{\sharp I-1}} \Lambda 
\ar[d]^\simeq \ar[r] & \lexp I \Lambda_I \ar[d]^\simeq \\
\Lambda^{\binom{\sharp I-1}{k-1}} \ar@{^{(}->}[r] &
\Lambda^{\binom{\sharp I}{k}} \ar[r] & \Lambda^{\binom{\sharp I}{k+1}} \ar[r] & \cdots \ar[r] &
\Lambda^{\binom{\sharp I}{1}} \ar[r] & \Lambda.
}$$
\end{prop}

\begin{proof}
As $\gr^k_!(\Psi(\Lambda)) (\frac{k-1}{2}) \hookrightarrow \gr^1_!(\Psi(\Lambda))$, we just need to
look at $\gr^1_!(\Psi(\Lambda))$. Recall that for every non empty $I \subset \{1,\cdots,r\}$, we have
$j_{I,!} \Lambda_I(\frac{\sharp I -1}{2}) \hookrightarrow \gr^1_!(\Psi(\Lambda))$ compatible with
their filtration $\Fill^\bullet$ and so a morphism of spectral sequences computing the $h^i z_I^*$ 
giving rise to the following commutative diagram
$$\xymatrix{
h^{\sharp I-r} z_I^* \lexp p j_{I,!*} \Lambda_I(\frac{\sharp I-1}{2}) \Bigr )
\ar[r] \ar[d] & 
h^{\sharp I+1-r} z_I^* \Bigl (
\bigoplus_{\myatop{J \supset I}{\sharp J \setminus I=1}} \lexp p j_{J,!*} \Lambda_J(\frac{\sharp I}{2}) \Bigr ) \ar[d] \\
h^{\sharp I-r} z_I^* \Bigl (
\bigoplus_{J:~\sharp J=\sharp I} \lexp p j_{J,!*} \Lambda_J(\frac{\sharp I-1}{2}) \Bigr )
\ar[r] & 
h^{\sharp I-r+1} z_I^* \Bigl (
\bigoplus_{J:~\sharp J=\sharp I+1} \lexp p j_{J,!*} \Lambda_J(\frac{\sharp I}{2}) \Bigr )
}$$
that is with the previous notation
$$\xymatrix{
\lexp I \Lambda \ar[r]^{f_I} \ar[d] & \bigoplus_{\myatop{J \supset I}{\sharp J \setminus I=1}} \lexp J \Lambda \ar[d] \\
\bigoplus_{J:~\sharp J=\sharp I} \lexp J \Lambda \ar[r]^{g_I} & 
\bigoplus_{J:~\sharp J=\sharp I+1} \lexp J \Lambda
}$$
where, up to modifying the isomorphisms $\lexp J \Lambda \simeq \Lambda$, as
$h^i z_I^* j_{I,!} \Lambda_I=(0)$, we have 
$f_I(1)=(1,\cdots,1)$. As moreover, for any $J$ such that $I \not \subset J$, every extension between
$\lexp p j_{I,!*} \Lambda_I(\frac{\sharp I-1}{2})$ and $\lexp p j_{J,!*} \Lambda_J(\frac{\sharp J-1}{2})$
is split, we then deduce that the image of $1 \in \lexp I \Lambda$ by $g_I$ is equal, up
to modify the isomorphisms $\lexp J \Lambda \simeq \Lambda$, to $(1,\cdots,1,0,\cdots,0)$
where we have $1$ (resp. $0$) for every $J \supset I$ (resp. $I \not \subset J$).

\end{proof}

Dually the cofiltration of stratification
$$\Psi(\Lambda)=\coFil_{*,r}(\Psi(\Lambda)) \twoheadrightarrow \coFil_{*,r-1}(\Psi(\Lambda))
\twoheadrightarrow \cdots \twoheadrightarrow \coFil_{*,0}(\Psi(\Lambda))=0,$$
with graded parts $\cogr_{*,k}(\Psi(\Lambda))$ having a filtration
\begin{multline*}
0=\Fill^{r-k}(\cogr_{*,k}(\Psi(\Lambda))) \hookrightarrow
\Fill^{r-k+1}(\cogr_{*,k}(\Psi(\Lambda))) \hookrightarrow \cdots  \\ 
\hookrightarrow \Fill^r(\cogr_{*,k}(\Psi(\Lambda)))=\cogr_{*,k}(\Psi(\Lambda))
\end{multline*}
with graded parts $\grr^h(\cogr_{*,k}(\Psi(\Lambda))) \simeq \lexp p j^{(h)}_{!*} \Lambda^{(h)}
(\frac{1-h+2(r-k))}{2}).$
Moreover the adjunction morphisms give monomorphism with torsion free cokernels
$$\cogr_{k,*}(\Psi(\Lambda)) \hookrightarrow j^{(r-k+1)}_! \Lambda^{(r-k+1)} (\frac{r-k}{2}).$$
As above, using (\ref{eq-jFill2}), the extensions inside $\Fill^\bullet(\cogr_{*,k}(\Psi(\Lambda)))$ 
do not split so that we have monomorphisms with torsion free cokernels
\begin{equation} \label{eq-psiFill2}
\cogr_{*,h}(\Psi(\Lambda))/\Fill^k(\cogr_{*,h}(\Psi(\Lambda))) \hookrightarrow
i^{(k+1)}_* j^{(k+1)}_* \Lambda^{(k+1)}(\frac{1-k+2(h-1)}{2}).
\end{equation}

Note finally that we never need to use the saturation process of definition \ref{defi-saturation}, for
$\Lambda=\overline \Zm_l$, when constructing the filtrations of stratification.

\begin{figure}[ht]
\centering
\input{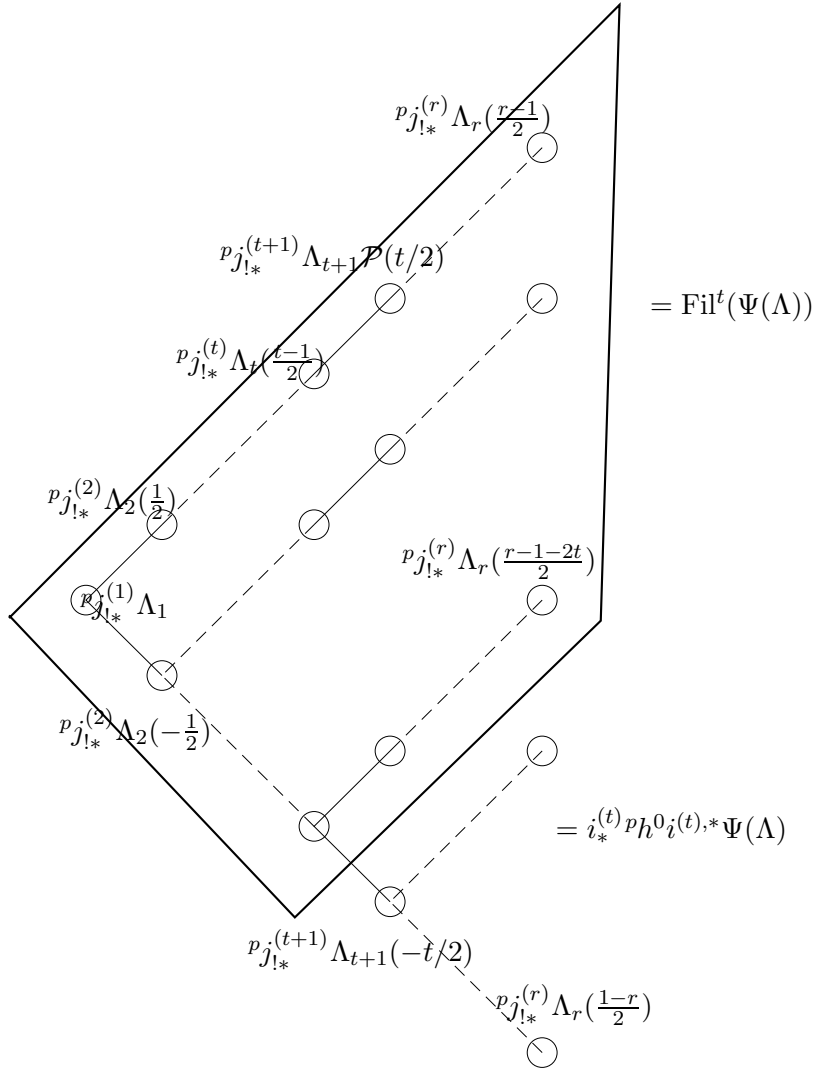}
\caption{\label{figure3} Filtration of stratification of $\Psi(\Lambda)$.}
\end{figure}

\begin{figure}[ht]
\centering
\input{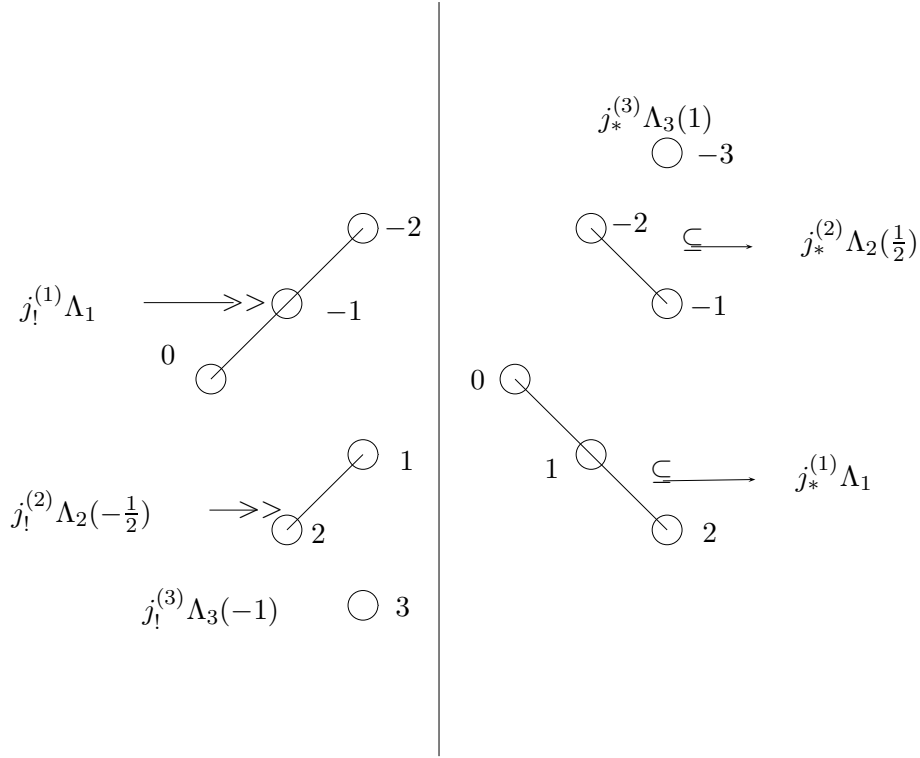}
\caption{\label{fig-filtration} Filtration (fig. on the left)
and cofiltration (fig. on the right) of stratification of $\Psi(\Lambda)$ with $r=3$: the number
correspond to the weights.}
\end{figure}

Combining lemma \ref{lem-analogpsi} with proposition \ref{prop-ij} we then deduce the following result
which was proved in \cite{dat-lem} theorem 2.2 or \cite{zheng} lemma 5.6 using different tools.

\begin{prop} \label{prop-datlem}
(cf. \cite{dat-lem} theorem 2.2 or \cite{zheng} lemma 5.6) \\
For every $I \subset \{ 1,\cdots,r \}$, the adjunction morphism
$$i_{I,*} i_I^* \Psi(\Lambda) \longrightarrow j_{I,*} j_I^* \Psi(\Lambda)$$
is an isomorphism.
\end{prop}

\rem $j_I$ being affine so that $j_{I,*} j_I^* \Psi(\Lambda)$ is perverse, the map of the
proposition factorizes through $\lexp p h^0 i_I^* \Psi(\Lambda)$.

\rem as
$$i_*^{(t)}\lexp p h^{0} i^{(t),*} \Psi(\Lambda) \simeq \Psi(\Lambda)/\Fil^t(\Psi(\Lambda))$$
we can obtain an explicit filtration of it with graded parts the perverse sheaves
$\lexp p j^{(k)}_{!*} \Lambda^{(k)}$.
Similarly
$$i_*^{(t)}\lexp p h^{-1} i^{(t),*} \Psi(\Lambda) \simeq \ker_\FC \bigl ( 
j^{(t)}_! j^{(t),*} \Psi(\Lambda) \twoheadrightarrow \Fil^t(\Psi(\Lambda)) \bigr ),$$ 
can also be described explicitly. 
Dually we can also give a filtration of
 $$i_*^{(t)}\lexp p h^{1} i^{(t),!} \Psi(\Lambda)=j_*^{(t)} j^{(t),*} \Psi(\Lambda)/\coFil_{*,t}(\Psi(\Lambda))$$
 and
 $$i_*^{(t)}\lexp p h^{0} i^{(t),!} \Psi(\Lambda)=\ker_\FC(\Psi(\Lambda) \longrightarrow  \coFil_{*,t}
 (\Psi(\Lambda))).$$

\section{The nilpotent monodromy operator}
\label{para-N}

First start with the knowledge of the existence of a non zero nilpotent monodromy operator over
$\Lambda=\overline \Qm_l$
$$N: \Psi(\Lambda) \longrightarrow \Psi(\Lambda)(1).$$
Note that the kernel of $N$ contains necessary $\Fil^1_!(\Psi(\Lambda))$: indeed for 
any of its graded parts
$\grr^k(\gr^1_!(\Psi(\Lambda)))$ then $\grr^k(\gr^1_!(\Psi(\Lambda)))(1)$ is not a constituant of
$\Psi(\Lambda)$.

As the socle of $\Psi(\Lambda)$
is $i^{(r)}_* \lexp p j^{(r)}_{!*} \Lambda^{(r)}(\frac{r-1}{2})$, then
$i^{(r)}_* \lexp p j^{(r)}_{!*} \Lambda^{(r)}(\frac{r-3}{2})$ does not belong to the kernel of $N$: note that
it appears with multiplicity one in $\Psi(\Lambda)$ as all the other constituants.

\begin{lem}
The socle of $i_*^{(1)}\lexp p h^{0} i^{(1),*} \Psi(\Lambda)$ is equal to 
$i^{(r)}_* \lexp p j^{(r)}_{!*} \Lambda^{(r)}(\frac{r-3}{2})$.
\end{lem}

\begin{proof}
We have already seen that the socle of the $\grr^h(\gr^k_!(\Psi(\Lambda)))$ are equal to
$i^{(r)}_* \lexp p j^{(r)}_{!*} \Lambda^{(r)}(\frac{r-1+2(k-1)}{2})$. Dually we already remarked that the
socle of the 
$$\grr^h(\cogr_{*,k}(\Psi(\Lambda)))/ i^{(k)}_*\lexp p j^{(k)}_{!*} \Lambda^{(k)}
(\frac{k-1}{2}),$$ 
is $i^{(k+1)}_* \lexp p j^{(k+1)}_{!*} \Lambda^{(k+1)}(\frac{k-3}{2})$.
Then the only possibility for being a subspace of
$$i_*^{(1)}\lexp p h^{0} i^{(1),*} \Psi(\Lambda) \simeq \Psi(\Lambda)/\Fil^1_!(\Psi(\Lambda))$$
is $i^{(r)}_* \lexp p j^{(r)}_{!*} \Lambda^{(r)}(\frac{r-3}{2})$.
\end{proof}

Finally we then obtain that 
$$\ker N= \Fil^1_!(\Psi(\Lambda)),$$
and $N$ induces an isomorphism
$$\xymatrix{
N:i_*^{(1)}\lexp p h^{0} i^{(1),*} \Psi(\Lambda) \ar[r]^\sim \ar[d]^\sim & 
i_*^{(1)}\lexp p h^{0} i^{(1),!} \Psi(\Lambda)(1) \ar[d]^\sim \\
\Psi(\Lambda) \simeq \Psi(\Lambda)/\Fil^1_!(\Psi(\Lambda)) &
\ker(\Psi(\Lambda) \rightarrow \coFil_{r-1,*}(\Psi(\Lambda)))(1).
}$$

As the socle of $i_*^{(1)}\lexp p h^{0} i^{(1),*} \Psi(\overline \Qm_l)$ and
$i_*^{(1)}\lexp p h^{0} i^{(1),!} \Psi(\overline \Qm_l)(1)$ are irreducible
isomorphic to $i^{(r)}_* \lexp p j^{(r)}_{!*} \overline \Qm_l^{(r)}(\frac{r-1}{2})$,
this isomorphism is unique up to a scalar so that we can find an $\overline \Zm_l$
version of $N$ inducing a $\overline \Zm_l$-isomorphism on
$i^{(r)}_* \lexp p j^{(r)}_{!*} \overline \Zm_l^{(r)}(\frac{r-1}{2})$. By reducing modulo $l$,
we obtain a $\overline \Fm_l$-version 
$$\overline N: \Psi(\overline \Fm_l) \longrightarrow \Psi(\overline \Fm_l)(1)$$
factorizing through
\begin{equation} \label{eq-barN}
\overline N:i_*^{(1)}\lexp p h^{0} i^{(1),*} \Psi(\overline \Fm_l) \longrightarrow
i_*^{(1)}\lexp p h^{0} i^{(1),!} \Psi(\overline \Fm_l)(1).
\end{equation}
As the socle of $i_*^{(1)}\lexp p h^{0} i^{(1),*} \Psi(\overline \Fm_l)$ and
$i_*^{(1)}\lexp p h^{0} i^{(1),!} \Psi(\overline \Fm_l)(1)$ are irreducible both isomorphic
to $i^{(r)}_* \lexp p j^{(r)}_{!*} \overline \Fm_l^{(r)}(\frac{r-1}{2})$ and as we choose $\overline N$ so that
it induces an isomorphism on this socle, then (\ref{eq-barN}) is an isomorphism.

Finally we obtain version of $N$ over $\Lambda \in \{ \overline \Qm_l,\overline \Zm_l,\overline \Fm_l \}$
which are nilpotent with order of nilpotency equal to $r$.

\rem It is also possible to prove directly that $i_*^{(1)}\lexp p h^{0} i^{(1),*} \Psi(\Lambda)$
and $i_*^{(1)}\lexp p h^{0} i^{(1),!} \Psi(\Lambda)(1)$ are isomorphic using the filtration of
stratification. One can argue by induction starting from the fact that
$\Fil^2_!$ of these two perverse sheaves are isomorphic obtained as a quotient of
$j^{(2)}_! \Lambda^{(2)}(\frac{1}{2})$. Then one can define $N$ directly
$$\xymatrix{
N: \Psi(\Lambda) \ar@{->>}[r] & \Psi(\Lambda)/\Fil^1_!(\Psi(\Lambda) \simeq
i_*^{(1)}\lexp p h^{0} i^{(1),*} \Psi(\Lambda) \ar[dl]^\sim \\
i_*^{(1)}\lexp p h^{0} i^{(1),!} \Psi(\Lambda)(1) \ar[r]^\sim &
\ker \Bigl ( \Psi(\Lambda) \rightarrow \cogr_{*,r-1}(\Psi(\Lambda)) \Bigr )(1) \ar@{^{(}->}[d] \\
& \Psi(\Lambda)(1)
}$$

\section{Cancellation of the $\Fm_l$-cohomology}
\label{para-application}

Having in mind application to the cohomology of Shimura varieties with coefficient in $\overline \Fm_l$,
we give in this section a baby example to motivate the conjecture 6.15 of \cite{K-S}, that we expect
to be used in more general situations.

To illustrate the idea, we consider the geometric setting
developed in \cite{Haines2012} \S 3.4 where we have a proper morphism
$$\pi: \AC_1 \longrightarrow \AC_0$$
with $\AC_0$ having strict semi stable reduction over $\spec \Zm_p$, i.e. its geometric special fiber
being as before a divisor with simple normal crossings
$$\AC_{0,\bar s} = \bigcup_{i=1}^r \AC_{0,i}.$$
As in the previous sections for
any non empty subset $I \subset \{ 1,\cdots,r \}$, put $\AC_{0,I}=\bigcup_{i \in I} \AC_{0,i}$.
Concerning $\AC_1$,
\begin{itemize}
\item $\pi$ is a connected Galois covering in the generic fiber with group of deck transformations 
$T(\Fm_p)$ where $T=\prod_{i=1}^r \Gm_m$; 

\item the special fiber $\AC_1 \otimes_{\Zm_p} \Fm_p$ is a divisor with normal crossing 
with all branches of multiplicity $p-1$;

\item for every non empty $I \subset \{ 1,\cdots,r \}$ let $\AC_{1,I}:=\pi^{-1}(\AC_{0,I})$. The restriction 
$\pi_I$ of $\pi$ to $\AC_{1,I}$ is an \'etale Galois covering over $\AC_{0,I}$ with group of deck
transformations $T^I(\Fm_p)$ where $T^I=\prod_{i \not \in S} \Gm_m$. 
\end{itemize}
For $x \in \AC_{0,I}(\overline \Fm_p)$, the covering $\pi:\AC_1 \longrightarrow \AC_0$ is an \'etale
local neighborhood of $x$ isomorphic to the morphism
$$\spec \Zm_p^{nr} [T'_1,\cdots,T'_r]/\Bigl ( (\prod_{i \in S} T'_i)^{p-1}- \varpi_p \Bigr ) \longrightarrow
\spec \Zm_p^{nr} [T_1,\cdots,T_r]/(\prod_{i \in S} T_i- \varpi_p )$$
defined by $T_i \mapsto (T'_i)^{p-1}$. Moreover we have a canonical surjective morphism
$$\pi_1(\AC_{0,I},x) \longrightarrow \Tm^I(\Fm_p).$$

If $\Lambda$ is large enough,
we have a direct sum decomposition into eigenspaces for the $T(\Fm_p)$ action
$$\pi_{\eta,*}(\Lambda)=\bigoplus_{\chi:T(\Fm_p) \rightarrow \Lambda_\times} \Lambda_\chi,$$
where $\Lambda_\chi$ is a local system of rank one on $\AC_{0,\eta}$ and if $\Psi_1(\Lambda)$
denotes the complex of nearby cycles on $\AC_{1} \otimes_{\Zm_p} \overline \Fm_p$, 
we have a direct sum decomposition
$$\pi_* \Psi_1(\Lambda)=\bigoplus_{\chi:T(\Fm_p) \rightarrow \Lambda_\times} \Psi_\chi(\Lambda),$$
with $\Psi_\chi(\Lambda):=\Psi(\Lambda_\chi)$.
In \cite{Haines2012} (5.4.6) the authors gives the following description of $\Psi_\chi(\Lambda)$:
\begin{equation} \label{eq-iso-I}
i_I^* \Psi_\chi(\Lambda) \simeq (\LC_I \otimes \Lambda_{\chi^I}) \otimes i_I^* \Psi(\Lambda)
\end{equation}
when $\chi$ is induced by $\chi^I:T^I(\Fm_p) \longrightarrow \Lambda^\times$ and is zero
otherwise.

In \cite{K-S} conjecture 6.15, for a reductive group $G/\Qm$ giving rise to a Shimura variety 
$\sh_K \rightarrow \spec F$ of level
$K$, an open compact subgroup of $G(\Am_\Qm)$, over its reflex field and of relative dimension
$d-1$, the authors try to predict what happen for its
$\overline \Fm_l$-cohomology groups. 
More precisely they look at the action of $G(\Qm_p)$ on 
$$\lim_{\myatop{\rightarrow}{K_p}} H^i(\sh_{K^pK_p},\overline \Fm_l[d-1]),$$
where $K_p$ describe the set of open compact subgroups of $G(\Qm_p)$ and propose
some constraints about its Arthur local parameter. Let $\Gamma_1(p)$ be a pro-$p$-Iwahori subgroup
i.e. the pro-$p$-unipotent radical of a Iwahori subgroup $\Gamma_0(p)$.
As taking invariants under $\Gamma_1(p)$
(resp. $\Gamma_0(p)$) is an exact
functor, representations of $G(\Qm_p)$ having non trivial invariant vectors under $\Gamma_1(p)$
can be detected when looking at level $\Gamma_1(p)$.

We resume the geometric situation of \cite{Haines2012} where
the algebraic group $G$ over $\Qm$ with values $\Qm$-algebras is defined as
$$G(R)=\{ x \in (D \otimes_\Qm R)^\times : x.x^* \in R^\times \},$$
where $D$ is a central simple $F$-algebra for some imaginary quadratic extension $F/\Qm$ and
$*$ is an involution of $D$ inducing the non trivial automorphism of $F$. We suppose that
$p$ splits in $F$ and 
$$D \otimes_\Qm \Qm_p \simeq \Mm_d(\Qm_p) \times \Mm_d(\Qm_l)^{op},$$
where $*$ transports to $(X,Y) \mapsto (Y^t,X^t)$. Note that
$$G \times_{\spec \Qm} \spec \Qm_p \simeq GL_d \times \Gm_m.$$
Then the Shimura varieties with respective level $K^p \Gamma_1(p)$ and 
$K^p \Gamma_0(p)$, verify the same hypothesis as $\pi:\AC_1 \longrightarrow \AC_0$ and we keep
the same notation $\AC_1$ and $\AC_0$ for them. 

\noindent \textit{Affiness hypothesis}: We moreover suppose that, for every
non empty subset $I \subset \{1,\cdots,r \}$, the $\overline \Fm_l$-cohomology groups of
$\AC_{0,I}^0$ are concentrated between $0$ and its dimension (as if they were affine).
 
This is the main restrictive hypothesis however there are some results in this direction 
which could be used.
\begin{itemize}
\item The hypothesis is verified for the Shimura varieties of Harris-Taylor type either in the split
and the inert cases.

\item According to \cite{mao} corollary 5.3, if the Shimura variety is compact then the central
leaves are affine. 

\item If moreover we are in the fully Hodge-Newton decomposable case,
then by \cite{GHN} theorem 7.4 the non basic Newton stratum are a finite union of central leaves and
are so also affine.
Loc. cit. also tells us that the EKOR stratification is a refinement of the
Newton stratification: note also that in the Iwahori level the EKOR stratification coincides with the
KR stratification. We then obtain the affiness of the non basic $Y_I^0$.

\item Another way to obtain such cancellation of the cohomology, notably for the basic locus, 
is to pass through the local model
which could be related to parabolic Deligne Lusztig varieties which are known to verify the previous hypothesis
\end{itemize}

\begin{lem} Let $\chi: \Tm(\Fm_p) \longrightarrow \overline \Fm_l^\times$ which factorizes through
$\chi^I:\Tm^I(\Fm_p) \longrightarrow \overline \Fm_l^\times$ for some fixed non empty subset
$I \subset \{ 1,\cdots,r \}$. Then for every $J \subset \{ 1,\cdots,r \}$
with $J \supset I$, if $H^{i}(\AC_{0,J},\lexp p j_{J,!*} \overline \Fm_{l,\chi^I})$ is non zero then
$|i| \leq \sharp I-\sharp J$.
\end{lem}

\begin{proof}
Note that  $H^{i}(\AC_{0,J},\lexp p j_{J,!*} \overline \Fm_{l,\chi^I})=(0)$ for every $i$ if 
$I \not \subset J$.
Using this property and the filtration of stratification of $j_{I,!} \overline \Fm_{l,\chi^I}$ of \S \ref{para-fil!}
we then deduce that 
$$H^i(\AC_0,j_{J,!} \overline \Fm_{l,\chi^I}) \simeq H^i(\AC_0^J, \lexp p j_{J,!*} \overline \Fm_{l,\chi^I}).$$
The affiness hypothesis implies it is zero if $i<0$ and the Grothendieck-Verdier duality on the left
hand side tells us then it is also zero for $i>0$.

We then argue by induction on $\sharp J$ from $\sharp I$ to $1$: we then suppose that for any
$J$ such that $h < \sharp J \leq \sharp I$, we have
$H^i(\AC_{0,\bar s},\lexp p j_{J,!*} \overline \Fm_{l,\chi^I})$ is zero if $|i| > \sharp I - \sharp J$ and
we consider $J_0$ with $\sharp J_0=h$. We then look at the spectral sequence computing
the $H^i(\AC_{0,\bar s},j_{J_0,!} \overline \Fm_{l,\chi^I})$ using the filtration of stratification
$\Fil^\bullet_!(j_{J_0,!} \overline \Fm_{l,\chi^I})$ of \S \ref{para-fil!} 
$$E_1^{i,j}=H^{i+j}(\AC_{0,\bar s},\gr^{-j}_!(j_{J_0,!} \overline \Fm_{l,\chi^I})) \Rightarrow 
H^{i+j}(\AC_{0,\bar s},j_{J_0,!} \overline \Fm_{l,\chi^I}),$$
as illustrated in the figure \ref{fig-ss-Hi} where we denote by 
$$j_{J_0}^{(k)}:\bigcup_{\myatop{J \supset J_0}{\sharp J=k}} \AC_{0,J}^0 \hookrightarrow 
\AC_{0,\bar s}.$$

\rem Using purity note that this spectral sequence degenerates at $E_2$.

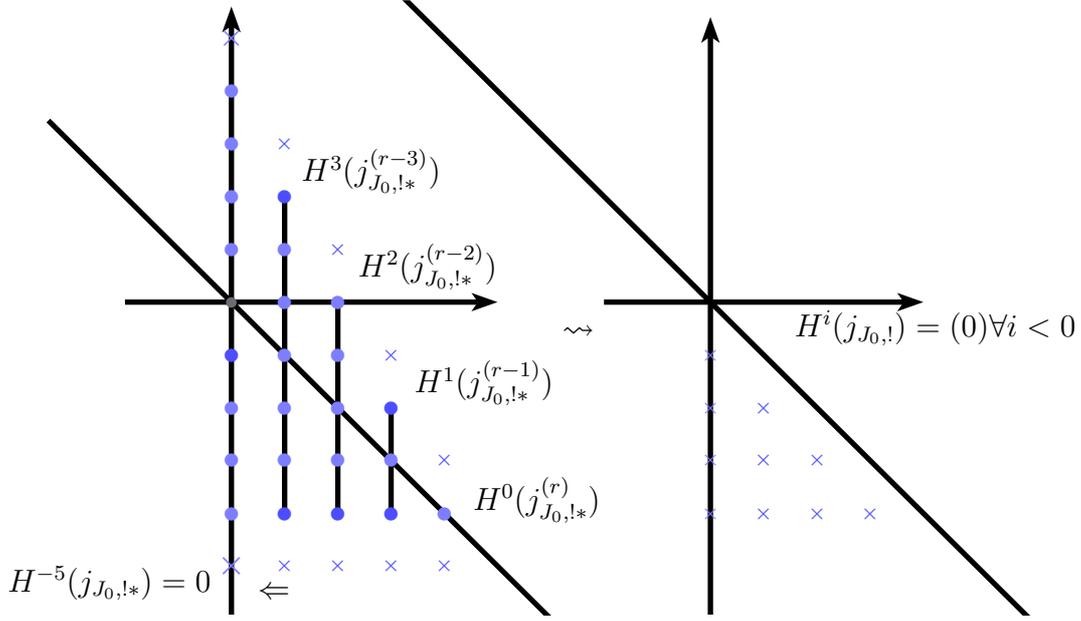
\begin{figure}[ht]
\centering
\input{fig-ss-Hi.tex}
\caption{\label{fig-ss-Hi} Spectral sequence computing the $H^i(j_{I,!} \Lambda_I)$: abutment is on the
right.}
\end{figure}
As $\AC_{0,J_0}^0$ is affine then $E_\oo^{i,j}=(0)$ for $i+j <0$ and as, by the induction hypothesis,
all the $E_1^{i,j}$ are zero if $i+j \leq h- \sharp I$ we obtain 
$H^{\sharp I - \sharp J_0-1}(\AC_{0,\bar s}, \lexp p j_{J_0,!*} \overline \Fm_{l,\chi^I})=(0)$.
\end{proof}

\begin{prop}
Let $i>0$ such $H^{-i}(\AC_{1,\bar s},\Psi_1(\overline \Fm_l)) \neq (0)$ then every irreducible subquotient
seen as a representation of $T(\Fm_p)$ factors through $T^I(\Fm_p)$ for some $I \subset \{1,\cdots,r \}$
with $\sharp I \geq i+1$.
\end{prop}

\begin{proof}
Consider a non empty subset $I \subset \{ 1,\cdots,r \}$ with $\sharp I \leq i$ and let 
$\chi:T^I(\Fm_p) \longrightarrow \overline \Fm_l^\times$. 
The filtration of stratification $\Fil^\bullet_!(\Psi(\overline \Fm_{l,\chi}))$ gives a spectral sequence
$$E_1^{p,q}=H^{p+q}(\AC_{0,\bar s},\gr^{-q}_!(\Psi((\overline \Fm_{l,\chi})))) \Rightarrow 
H^{p+q}(\AC_{0,\bar s},\Psi(\overline \Fm_{l,\chi})).$$
The filtration $\Fill^\bullet(\gr^k_!(\Psi(\overline \Fm_{l,\chi})))$ gives in the same way a spectral sequence
$$EE_1^{p,q}=H^{p+q}(\AC_{0,\bar s},\grr^{-q}(\gr^k_!(\Psi(\overline \Fm_{l,\chi})))) \Rightarrow 
H^{p+q}(\AC_{0,\bar s},\gr^{^k}_!(\Psi(\overline \Fm_{l,\chi}))).$$
Note that (\ref{eq-iso-I}) tells us that the cohomology of $\Psi_\chi(\overline \Fm_l)$ is supported on $\bigcup_{J \supset I} \AC_{0,J}^0$, so that, in the initial terms of this last spectral sequence, 
only appears the cohomology groups of the $\lexp p j_{J,!*} \overline \Fm_{l,\chi_I}$
for $J \supset I$. The previous lemma then tells us that all the $EE_1^{p,q}$ are zero for
$|p+q| > \sharp I \geq i$ and so it is the same for the $E_\oo^{p,q}$.

\end{proof}

\begin{defin}
A character $\chi$ is said generic if it does not factorise through any $T^I(\Fm_p)$ for
$\sharp I \geq 2$. 
\end{defin}

\begin{prop} For $\chi$ generic, the $\chi$ part of the $\overline \Fm_l$-cohomology of 
$\AC_1$ is concentrated in middle degree.
\end{prop}

\rem This is essentially the approach followed in \cite{boyer-imj}.

\bibliography{Biblio}
\bibliographystyle{amsalpha}

\end{document}

%% file: fig-ss-N.tex
\newrgbcolor{ududff}{0.30196078431372547 0.30196078431372547 1}
\newrgbcolor{xdxdff}{0.49019607843137253 0.49019607843137253 1}
\psset{xunit=1cm,yunit=1cm,algebraic=true,dimen=middle,dotstyle=o,dotsize=5pt 0,linewidth=2pt,arrowsize=3pt 2,arrowinset=0.25}
\begin{pspicture*}(-8.6,-5.5)(1.2,5.5)
\pscircle[linewidth=2pt](-6,0){0.3}
\pscircle[linewidth=2pt](-5,1){0.3}
\pscircle[linewidth=2pt](-5,-1){0.3}
\pscircle[linewidth=2pt](-4,2){0.3}
\pscircle[linewidth=2pt](-4,0){0.3}
\pscircle[linewidth=2pt](-4,-2){0.3}
\pscircle[linewidth=2pt](-3,3){0.3}
\pscircle[linewidth=2pt](-3,1){0.3}
\pscircle[linewidth=2pt](-3,-1){0.3}
\pscircle[linewidth=2pt](-3,-3){0.3}
\psline[linewidth=2pt](-6,0)(-3,3)
\psline[linewidth=2pt](-5,-1)(-3,1)
\psline[linewidth=2pt](-4,-2)(-3,-1)
\psline[linewidth=2pt](-6,0)(-3,-3)
\psline[linewidth=2pt](-5,1)(-3,-1)
\psline[linewidth=2pt](-4,2)(-3,1)
\parametricplot[linewidth=2pt]{-1.57}{1.57}{1*1*cos(t)+0*1*sin(t)+-3|0*1*cos(t)+1*1*sin(t)+0}
\parametricplot[linewidth=2pt]{-1.57}{1.57}{1*1*cos(t)+0*1*sin(t)+-3|0*1*cos(t)+1*1*sin(t)+-2}
\parametricplot[linewidth=2pt]{-1.57}{1.57}{1*1*cos(t)+0*1*sin(t)+-3|0*1*cos(t)+1*1*sin(t)+2}
\parametricplot[linewidth=2pt]{-0.42}{0.42}{1*2.44*cos(t)+0*2.44*sin(t)+-7.23|0*2.44*cos(t)+1*2.44*sin(t)+-0.01}
\parametricplot[linewidth=2pt]{-0.49}{0.44}{1*2.23*cos(t)+0*2.23*sin(t)+-5.97|0*2.23*cos(t)+1*2.23*sin(t)+1.05}
\parametricplot[linewidth=2pt]{-0.566}{0.522}{1*1.955*cos(t)+0*1.955*sin(t)+-5.65|0*1.955*cos(t)+1*1.955*sin(t)+-0.95}
\rput[tl](-2.19,0.36){$\bigwedge$}
\rput[tl](-2.19,2.36){$\bigwedge$}
\rput[tl](-2.19,-1.64){$\bigwedge$}
\rput[tl](-4.97,0.36){$\bigwedge$}
\rput[tl](-3.9,1.36){$\bigwedge$}
\rput[tl](-3.85,-0.64){$\bigwedge$}

\rput[tl](-1.69,0.36){N}
\rput[tl](-1.74,2.62){N}
\rput[tl](-1.643,-1.86){N}
\rput[tl](-5.3,0.343){N}
\rput[tl](-4.211,1.246){N}
\rput[tl](-4.142,-0.749){N}
\rput[tl](-7.2,-0.8){$gr^1_!$}
\rput[tl](-6,-1.7){$gr^2_!$}
\rput[tl](-5.1,-2.7){$gr^3_!$}
\rput[tl](-7.38,1.2){$cogr_{*,4}$}
\rput[tl](-6.38,2){$cogr_{3,*}$}
\rput[tl](-5.38,3){$cogr_{*,2}$}
\psline[linewidth=2pt]{->}(-6.76,-0.71)(-6.34,-0.33)
\psline[linewidth=2pt]{->}(-5.66,-1.63)(-5.34,-1.33)
\psline[linewidth=2pt]{->}(-4.68,-2.53)(-4.255,-2.157)
\psline[linewidth=2pt]{->}(-6.66,0.71)(-6.28,0.31)
\psline[linewidth=2pt]{->}(-5.72,1.59)(-5.38,1.31)
\psline[linewidth=2pt]{->}(-4.8,2.47)(-4.4,2.17)
\rput[tl](-7.2,-2.345){$j^{(1)}_!$}
\rput[tl](-6,-3.4){$j^{(2)}_!$}
\rput[tl](-5,-4.4){$j^{(3)}_!$}
\rput[tl](-7.38,3){$j^{(1)}_*$}
\rput[tl](-6.38,3.7){$j^{(2)}_*$}
\rput[tl](-5.38,4.7){$j^{(3)}_*$}
\psline[linewidth=2pt]{->>}(-7,-2.1)(-7,-1.35)
\psline[linewidth=2pt]{->>}(-5.8,-3.2)(-5.8,-2.2)
\psline[linewidth=2pt]{->>}(-4.7,-4.2)(-4.7,-3.3)
\psline[linewidth=2pt]{->}(-6.98,1.3)(-6.98,2.18)
\psline[linewidth=2pt]{->}(-6.12,2.148)(-6.12,2.999)
\psline[linewidth=2pt]{->}(-5,3.2)(-5,3.9)

\rput[tl](-7,1.5){$\smallsmile$}
\rput[tl](-6.14,2.3){$\smallsmile$}
\rput[tl](-5,3.4){$\smallsmile$}

\begin{scriptsize}
\psdots[dotstyle=*,linecolor=ududff](-6,0)
\psdots[dotstyle=*,linecolor=ududff](-3,3)
\psdots[dotstyle=*,linecolor=ududff](-3,-3)
\psdots[dotstyle=*,linecolor=ududff](-5,1)
\psdots[dotstyle=*,linecolor=ududff](-5,-1)
\psdots[dotstyle=*,linecolor=ududff](-4,2)
\psdots[dotstyle=*,linecolor=xdxdff](-4,0)
\psdots[dotstyle=*,linecolor=ududff](-4,-2)
\psdots[dotstyle=*,linecolor=ududff](-3,3)
\psdots[dotstyle=*,linecolor=ududff](-3,1)
\psdots[dotstyle=*,linecolor=ududff](-3,-1)
\psdots[dotstyle=*,linecolor=ududff](-3,-3)
\end{scriptsize}
\end{pspicture*}

%% file: fig-hi-rPsi.tex
\newrgbcolor{ududff}{0.30196078431372547 0.30196078431372547 1}
\psset{xunit=1cm,yunit=1cm,algebraic=true,dimen=middle,dotstyle=o,dotsize=5pt 0,linewidth=2pt,arrowsize=3pt 2,arrowinset=0.25}
\begin{pspicture*}(-11.3,-3.41)(-0.456,4.55)
\psline[linewidth=2pt]{->}(-10,1)(-9,2)
\psline[linewidth=2pt]{->}(-9,2)(-8,3)
\psline[linewidth=2pt]{->}(-8,3)(-7,4)
\psline[linewidth=2pt]{->}(-9,0)(-8,1)
\psline[linewidth=2pt]{->}(-8,1)(-7,2)
\psline[linewidth=2pt]{->}(-8,-1)(-7,0)
\rput[tl](-11,1.5){$\Lambda^4(\frac{3}{2})$}
\rput[tl](-10,2.5){$\Lambda^6(\frac{3}{2})$}
\rput[tl](-9,3.5){$\Lambda^4(\frac{3}{2})$}
\rput[tl](-8,4.5){$\Lambda(\frac{3}{2})$}
\rput[tl](-10,0.5){$\Lambda^6(\frac{1}{2})$}
\rput[tl](-9,1.5){$\Lambda^4(\frac{1}{2})$}
\rput[tl](-8,2.5){$\Lambda(\frac{1}{2})$}
\rput[tl](-9.2,-1){$\Lambda^4(\frac{-1}{2})$}
\rput[tl](-8,0.5){$\Lambda(\frac{-1}{2})$}
\rput[tl](-6,1.8){$\leadsto$}
\rput[tl](-6,1){$\Lambda(\frac{3}{2})$}
\rput[tl](-5,0){$\Lambda^3(\frac{1}{2})$}
\rput[tl](-4.2,-1){$\Lambda^3(\frac{-1}{2})$}
\rput[tl](-3.2,-2){$\Lambda(\frac{-3}{2})$}
\rput[tl](-8.2,-2){$\Lambda(\frac{-3}{2})$}
\psline[linewidth=2pt,linestyle=dotted](-6,-3.4)(-6,4.5)
\begin{scriptsize}
\psdots[dotstyle=*,linecolor=ududff](-10,1)
\psdots[dotstyle=*,linecolor=ududff](-9,2)
\psdots[dotstyle=*,linecolor=ududff](-8,3)
\psdots[dotstyle=*,linecolor=ududff](-7,4)
\psdots[dotstyle=*,linecolor=ududff](-9,0)
\psdots[dotstyle=*,linecolor=ududff](-8,1)
\psdots[dotstyle=*,linecolor=ududff](-7,2)
\psdots[dotstyle=*,linecolor=ududff](-8,-1)
\psdots[dotstyle=*,linecolor=ududff](-7,0)
\psdots[dotstyle=*,linecolor=ududff](-7,-2)
\psdots[dotstyle=*,linecolor=ududff](-5,1)
\psdots[dotstyle=*,linecolor=ududff](-4,0)
\psdots[dotstyle=*,linecolor=ududff](-3,-1)
\psdots[dotstyle=*,linecolor=ududff](-2,-2)
\end{scriptsize}
\end{pspicture*}

%% file: fig-ss-Hi.tex
\newrgbcolor{bubtba}{0.7 0.7 0.73}
\newrgbcolor{ududff}{0.3 0.3 1}
\newrgbcolor{xdxdff}{0.49 0.49 1}
\newrgbcolor{wewdxt}{0.43 0.42 0.45}
\psset{xunit=.7cm,yunit=.7cm,algebraic=true,dimen=middle,dotstyle=o,dotsize=5pt 0,linewidth=2pt,arrowsize=3pt 2,arrowinset=0.25}
\begin{pspicture*}(-8.5,-2.92)(18,8.7)
\psline[linewidth=2pt,arrows=->](-4,-2.92)(-4,8.6)
\psline[linewidth=2pt,arrows=->](-6,3)(1,3)
\psline[linewidth=2pt,arrows=->](3,3)(9,3)
\psplot[linewidth=2pt]{-7.44}{11.92}{(-4-4*x)/4}
\psline[linewidth=2pt](-3,-1)(-3,5)
\psline[linewidth=2pt](-2,-1)(-2,3)
\psline[linewidth=2pt](-1,-1)(-1,1)
\rput[tl](0.54,-0.32){$H^0(j^{(r)}_{J_0,!*})$}
\rput[tl](-0.56,1.9){$H^{1}(j^{(r-1)}_{J_0,!*})$}
\rput[tl](-1.62,4.1){$H^2(j^{(r-2)}_{J_0,!*})$}
\rput[tl](-2.68,5.9){$H^3(j^{(r-3)}_{J_0,!*})$}
\rput[tl](-8.2,-2){$H^{-5}(j_{J_0,!*})=0$}
\rput[tl](-3.5,-2.4){$\Leftarrow$}
\rput[tl](2.22,2.54){$\leadsto$}
\psline[linewidth=2pt,arrows=->](5,-2.92)(5,8.4)
\psplot[linewidth=2pt]{-7.44}{11.92}{(--32-4*x)/4}
\rput[tl](6.58,2.84){$H^i(j_{J_0,!})=(0) \forall i<0$}
\begin{scriptsize}
\psdots[dotstyle=*,linecolor=ududff](-3,-1)
\psdots[dotstyle=*,linecolor=ududff](-3,5)
\psdots[dotstyle=*,linecolor=ududff](-2,-1)
\psdots[dotstyle=*,linecolor=xdxdff](-2,3)
\psdots[dotstyle=*,linecolor=ududff](-1,-1)
\psdots[dotstyle=*,linecolor=ududff](-1,1)
\psdots[dotstyle=*,linecolor=xdxdff](0,-1)
\psdots[dotstyle=x,linecolor=ududff](-3,-2)
\psdots[dotstyle=x,linecolor=ududff](-2,-2)
\psdots[dotstyle=x,linecolor=ududff](-1,-2)
\psdots[dotstyle=x,linecolor=ududff](0,-2)
\psdots[dotstyle=x,linecolor=ududff](-3,6)
\psdots[dotstyle=x,linecolor=ududff](-2,4)
\psdots[dotstyle=x,linecolor=ududff](-1,2)
\psdots[dotstyle=x,linecolor=ududff](0,0)

\psdots[dotstyle=*,linecolor=xdxdff](-3,0)
\psdots[dotstyle=*,linecolor=xdxdff](-2,0)
\psdots[dotstyle=*,linecolor=xdxdff](-1,0)
\psdots[dotstyle=*,linecolor=xdxdff](-3,1)
\psdots[dotstyle=*,linecolor=xdxdff](-2,1)
\psdots[dotstyle=*,linecolor=xdxdff](-3,2)
\psdots[dotstyle=*,linecolor=xdxdff](-2,2)
\psdots[dotstyle=*,linecolor=xdxdff](-3,3)
\psdots[dotstyle=*,linecolor=xdxdff](-3,4)

\psdots[dotsize=8pt 0,dotstyle=x,linecolor=xdxdff](-4,-2)
\psdots[dotstyle=*,linecolor=xdxdff](-4,-1)
\psdots[dotstyle=*,linecolor=xdxdff](-4,0)
\psdots[dotstyle=*,linecolor=xdxdff](-4,1)
\psdots[dotstyle=*,linecolor=ududff](-4,2)
\psdots[dotsize=4pt 0,dotstyle=*,linecolor=wewdxt](-4,3)
\psdots[dotstyle=*,linecolor=xdxdff](-4,4)
\psdots[dotstyle=*,linecolor=xdxdff](-4,5)
\psdots[dotstyle=*,linecolor=xdxdff](-4,6)
\psdots[dotstyle=*,linecolor=xdxdff](-4,7)
\psdots[dotsize=7pt 0,dotstyle=x,linecolor=xdxdff](-4,8)
\psdots[dotstyle=x,linecolor=xdxdff](5,2)
\psdots[dotstyle=x,linecolor=ududff](6,1)
\psdots[dotstyle=x,linecolor=ududff](7,0)
\psdots[dotstyle=x,linecolor=ududff](8,-1)
\psdots[dotstyle=x,linecolor=ududff](7,-1)
\psdots[dotstyle=x,linecolor=ududff](6,-1)
\psdots[dotstyle=x,linecolor=ududff](6,0)
\psdots[dotstyle=x,linecolor=xdxdff](5,1)
\psdots[dotstyle=x,linecolor=xdxdff](5,0)
\psdots[dotstyle=x,linecolor=xdxdff](5,-1)
\end{scriptsize}
\end{pspicture*}